\documentclass[12pt]{amsart}

\usepackage[margin=1.25in]{geometry}

\usepackage{amsmath}
\usepackage{amsthm}
\usepackage{amssymb}
\usepackage{caption}
\usepackage[curve]{xypic}
\usepackage{cite}
\usepackage{mathrsfs}
\usepackage{enumitem}
\usepackage{url}
\usepackage{color}
\usepackage[scr=boondoxo]{mathalpha}
\usepackage[usenames,dvipsnames]{xcolor}
\usepackage{graphicx}
\usepackage{tikz}
\usepackage{tikz-cd}
\usepackage{hyperref} 
\usepackage{verbatim}
\usepackage{multirow}
\usepackage{subcaption}
\usepackage{graphicx}
\usepackage{algorithm}
\usepackage{algorithmic}
\usetikzlibrary{decorations.markings,shapes.geometric}
\usetikzlibrary{decorations.pathreplacing}
 \usetikzlibrary{calc,decorations.markings,decorations.pathreplacing}

\newcommand{\dotpt}[2]{\fill (#1,#2) circle (1.7pt);}
\newcommand{\openpt}[2]{%
  \draw[fill=white] (#1,#2) circle (2.1pt);}
\newcommand{\squarept}[2]{%
  \draw[fill=white]
    (#1-.055,#2-.055) rectangle (#1+.055,#2+.055);}
\newcommand{\crosspt}[2]{%
  \draw (#1-.045,#2-.045)--(#1+.045,#2+.045)
        (#1-.045,#2+.045)--(#1+.045,#2-.045);}

\newcommand{\ticks}[4]{%
  \begin{scope}[shift={(#1,#2)},rotate=#3]
    \foreach \j in {1,...,#4}{
      \draw ({.065*(\j-(#4+1)/2)},-.055)
         -- ({.065*(\j-(#4+1)/2)}, .055);
    }
  \end{scope}}

\newcommand{\chevron}[3]{%
  \begin{scope}[shift={(#1,#2)},rotate=#3]
    \draw (-.065,-.065)--(0,0)--(-.065,.065);
  \end{scope}}

\newcommand{\sphere}[4]{%
  \begin{scope}[shift={(#1,#2)},xscale=#3,yscale=#4]
    \path[surface] (0,0) circle (1);
    \draw (-1,0)
      .. controls (-.55,-.42) and (.55,-.42) .. (1,0);
    \draw[core] (-1,0)
      .. controls (-.55,.38) and (.55,.38) .. (1,0);
  \end{scope}}

\newcommand{\handle}[2]{%
  \begin{scope}[shift={(#1,#2)}]
    \draw (-.16,.06)
      .. controls (-.08,-.10) and (.08,-.10) .. (.16,.06);
    \draw (-.105,.005)
      .. controls (-.12,.15) and (.12,.15) .. (.105,.005);
  \end{scope}}

\setlist{itemsep=4pt, topsep=4pt}

\makeatletter

\def\chaptermark#1{}%whatever

\def\chapter{%
\if@openright\cleardoublepage\else\clearpage\fi
\thispagestyle{plain}\global\@topnum\z@
\@afterindenttrue \secdef\@chapter\@schapter}

\def\@chapter[#1]#2{\refstepcounter{chapter}%
\ifnum\c@secnumdepth<\z@ \let\@secnumber\@empty
\else \let\@secnumber\thechapter \fi
\typeout{\chaptername\space\@secnumber}%
\def\@toclevel{0}%
\ifx\chaptername\appendixname \@tocwriteb\tocappendix{chapter}{#2}%
\else \@tocwriteb\tocchapter{chapter}{#2}\fi
\chaptermark{#1}%
\addtocontents{lof}{\protect\addvspace{10\p@}}%
\addtocontents{lot}{\protect\addvspace{10\p@}}%
\@makechapterhead{#2}\@afterheading}
\def\@schapter#1{\typeout{#1}%
\let\@secnumber\@empty
\def\@toclevel{0}%
\ifx\chaptername\appendixname \@tocwriteb\tocappendix{chapter}{#1}%
\else \@tocwriteb\tocchapter{chapter}{#1}\fi
\chaptermark{#1}%
\addtocontents{lof}{\protect\addvspace{10\p@}}%
\addtocontents{lot}{\protect\addvspace{10\p@}}%
\@makeschapterhead{#1}\@afterheading}
\newcommand\chaptername{Chapter}

\def\@makechapterhead#1{\global\topskip 7.5pc\relax
\begingroup
\fontsize{\@xivpt}{18}\bfseries\centering
\ifnum\c@secnumdepth>\m@ne
  \leavevmode \hskip-\leftskip
  \rlap{\vbox to\z@{\vss
      \centerline{\normalsize\mdseries
          \uppercase\@xp{\chaptername}\enspace\thechapter}
      \vskip 3pc}}\hskip\leftskip\fi
 #1\par \endgroup
\skip@34\p@ \advance\skip@-\normalbaselineskip
\vskip\skip@ }
\def\@makeschapterhead#1{\global\topskip 7.5pc\relax
\begingroup
\fontsize{\@xivpt}{18}\bfseries\centering
#1\par \endgroup
\skip@34\p@ \advance\skip@-\normalbaselineskip
\vskip\skip@ }

\newcounter{chapter}

\newif\if@openright

\makeatother

\makeatletter % makes "@" a normal character so that it can be used in macros
\def\@cite#1#2{{\m@th\upshape\bfseries%
[{#1\if@tempswa{\m@th\upshape\mdseries, #2}\fi}]}}
\makeatother % returns "@" to its normal state

\theoremstyle{plain}
\newtheorem{thm}{Theorem}[section]% subsection
\newtheorem{cor}[thm]{Corollary}

\newtheorem{prop}[thm]{Proposition}
\newtheorem{lem}[thm]{Lemma}

\theoremstyle{definition}

\theoremstyle{remark}
\newtheorem{rem}[thm]{Remark}

\numberwithin{equation}{subsection}
\newcommand{\nc}{\newcommand}
\newcommand{\rnc}{\renewcommand}
\newcommand{\e}{\varepsilon}

\newcommand{\ann}[1]{\marginpar{\tiny{#1}}}
\nc\bA{\mathbb{A}}
\nc\bB{\mathbb{B}}
\nc\bC{\mathbb{C}}
\nc\bD{\mathbb{D}}
\nc\bE{\mathbb{E}}
\nc\bF{\mathbb{F}}
\nc\bG{\mathbb{G}}
\nc\bH{\mathbb{H}}
\nc\bI{\mathbb{I}}
\nc{\bJ}{\mathbb{J}} 
\nc\bK{\mathbb{K}}
\nc\bL{\mathbb{L}}
\nc\bM{\mathbb{M}}
\nc\bN{\mathbb{N}}
\nc\bO{\mathbb{O}}
\nc\bP{\mathbb{P}}
\nc\bQ{\mathbb{Q}}
\nc\bR{\mathbb{R}}
\nc\bS{\mathbb{S}}
\nc\bT{\mathbb{T}}
\nc\bU{\mathbb{U}}
\nc\bV{\mathbb{V}}
\nc\bW{\mathbb{W}}
\nc\bY{\mathbb{Y}}
\nc\bX{\mathbb{X}}
\nc\bZ{\mathbb{Z}}
\nc\cA{\mathcal{A}}
\nc\cB{\mathcal{B}}
\nc\cC{\mathcal{C}}
\rnc\cD{\mathcal{D}}
\nc\cE{\mathcal{E}}
\nc\cF{\mathcal{F}}
\nc\cG{\mathcal{G}}
\rnc\cH{\mathcal{H}}
\nc\cI{\mathcal{I}}
\nc{\cJ}{\mathcal{J}} 
\nc\cK{\mathcal{K}}
\rnc\cL{\mathcal{L}}
\nc\cM{\mathcal{M}}
\nc\cN{\mathcal{N}}
\nc\cO{\mathcal{O}}
\nc\cP{\mathcal{P}}
\nc\cQ{\mathcal{Q}}
\rnc\cR{\mathcal{R}}
\nc\cS{\mathcal{S}}
\nc\cT{\mathcal{T}}
\nc\cU{\mathcal{U}}
\nc\cV{\mathcal{V}}
\nc\cW{\mathcal{W}}
\nc\cY{\mathcal{Y}}
\nc\cX{\mathcal{X}}
\nc\cZ{\mathcal{Z}}
\nc\bfA{\mathbf{A}}
\nc\bfB{\mathbf{B}}
\nc\bfC{\mathbf{C}}
\nc\bfD{\mathbf{D}}
\nc\bfE{\mathbf{E}}
\nc\bfF{\mathbf{F}}
\nc\bfG{\mathbf{G}}
\nc\bfH{\mathbf{H}}
\nc\bfI{\mathbf{I}}
\nc{\bfJ}{\mathbf{J}} 
\nc\bfK{\mathbf{K}}
\nc\bfL{\mathbf{L}}
\nc\bfM{\mathbf{M}}
\nc\bfN{\mathbf{N}}
\nc\bfO{\mathbf{O}}
\nc\bfP{\mathbf{P}}
\nc\bfQ{\mathbf{Q}}
\nc\bfR{\mathbf{R}}
\nc\bfS{\mathbf{S}}
\nc\bfT{\mathbf{T}}
\nc\bfU{\mathbf{U}}
\nc\bfV{\mathbf{V}}
\nc\bfW{\mathbf{W}}
\nc\bfY{\mathbf{Y}}
\nc\bfX{\mathbf{X}}
\nc\bfZ{\mathbf{Z}}

\newcommand{\del}{\partial}

\nc{\dmo}{\DeclareMathOperator}
\nc{\wt}{\widetilde}
\rnc{\Re}{\operatorname{Re}}
\rnc{\Im}{\operatorname{Im}}
\rnc{\span}{\operatorname{span}}
\dmo{\rank}{rank}
\dmo{\End}{End}
\dmo{\Hom}{Hom}
\dmo{\Jac}{Jac}
\dmo{\Id}{Id}
\dmo{\Ann}{Ann}
\dmo{\Area}{Area}
\dmo{\CP}{\bC P^1}
\dmo{\rk}{rk}
\dmo{\rel}{rel}
\dmo{\ra}{\rightarrow}
\dmo{\AGL}{\mathrm{AGL}}
\dmo{\AO}{\mathrm{AO}}
\dmo{\Sym}{\mathrm{Sym}}
\dmo{\Hur}{\mathrm{Hur}}
\dmo{\Aut}{\mathrm{Aut}}
\dmo{\Mod}{\mathrm{Mod}}

\rnc{\Col}{\operatorname{Col}}

\nc{\ColOne}{\Col_{\bfC_1}}
\nc{\ColOneX}{\ColOne(X,\omega)}

\nc{\ColTwo}{\Col_{\bfC_2}}
\nc{\ColTwoX}{\ColTwo(X,\omega)}

\nc{\ColThree}{\Col_{\bfC_3}}
\nc{\ColThreeX}{\ColThree(X,\omega)}

\nc{\ColOneTwo}{\Col_{\bfC_1, \bfC_2}}
\nc{\ColOneTwoX}{\ColOneTwo(X,\omega)}

\nc{\ColOneThree}{\Col_{\bfC_1, \bfC_3}}
\nc{\ColOneThreeX}{\ColOneThree(X,\omega)}

\nc{\MOne}{\cM_{\bfC_1}}
\nc{\MTwo}{\cM_{\bfC_2}}
\nc{\MOneTwo}{\cM_{\bfC_1, \bfC_2}}

\nc{\MThree}{\cM_{\bfC_3}}

\nc{\MOneThree}{\cM_{\bfC_1, \bfC_3}}

\dmo{\For}{\cF}

\nc{\GL}{\mathrm{GL}^+(2, \bR)}

\newcommand{\markededge}[4]{%
  \draw[
    decoration={markings,
      mark=at position .62 with {%
        \ifnum#4=1
          \arrow{Stealth[length=3pt,width=3pt]}
        \fi
      },
      mark=at position .30 with {%
        \ifnum#3>0
          \foreach \k in {1,...,#3}{%
            \draw[line width=.5pt]
              ({(\k-(#3+1)/2)*1.9pt},-2pt)
              -- ({(\k-(#3+1)/2)*1.9pt},2pt);
          }%
        \fi
      }
    },
    postaction={decorate}
  ] (#1)--(#2);
}

\newcommand{\glueedge}[3]{\markededge{#1}{#2}{#3}{1}}
\newcommand{\sideedge}[3]{\markededge{#1}{#2}{#3}{0}}

\newcommand{\foldedge}[4]{%
  \coordinate (#3) at ($(#1)!.5!(#2)$);
  \glueedge{#1}{#3}{#4}
  \glueedge{#2}{#3}{#4}
}

\newcommand{\pole}[1]{%
  \node[inner sep=0pt,font=\scriptsize] at (#1) {$\times$};
}

\newcommand{\edgepole}[3]{%
  \path
    let \p1=($(#3)-(#2)$),
        \n1={atan2(\y1,\x1)}
    in node[
      rotate=\n1,
      inner sep=0pt,
      font=\scriptsize
    ] at (#1) {$\times$};
}

\newcommand{\zeropoint}[1]{%
  \fill (#1) circle (1.8pt);
}

\title[The converse of Kra's theorem]{The converse of Kra's Theorem}
\author[P. Apisa]{Paul Apisa}
\author[J. Aygun]{Juliet Aygun}
\author[D. Novikov]{Dmitry Novikov}
\author[A. Zimmer]{Andrew Zimmer}
\begin{document}

\begin{abstract}
A conjectured converse to a theorem of Kra states that a Teichm\"uller disk is a holomorphic retract of Teichm\"uller space if and only if it is generated by a quadratic differential whose zeros are even order when not at punctures. We prove this for $\cT_{0,n}$ and for $\cT_{1,n}$ for all $n$ with the only possible exceptions coming from the ornithorynque Teichm\"uller curve. In general, we show that the conjecture holds on a full measure set of all strata of quadratic differentials with five possible exceptions.  
\end{abstract}

\maketitle

%%%%%%%%%%%%%%%%%%% 
% TABLE OF CONTENTS
%%%%%%%%%%%%%%%%%%%
% allows subsections (depth 1) to be displayed in table of contents
\setcounter{tocdepth}{1} 
\tableofcontents

\vspace{-5mm}

\section{Introduction}\label{S:Intro}

Let $\mathcal T_{g,n}$ be the Teichm\"uller space of marked Riemann surfaces of genus $g$ with $n$ punctures. A \emph{Teichm\"uller disk}, also called a \emph{complex geodesic}, in $\mathcal T_{g,n}$ is a holomorphic embedding $\tau : \mathbb{H} \to \mathcal T_{g,n}$ which is an isometry with respect to the Poincar\'e metric on $\mathbb{H}$ and Teichm\"uller metric on $\mathcal T_{g,n}$. Say that the Teichm\"uller disk \emph{is a holomorphic retract} if there is a holomorphic map $r\colon\mathcal{T}_{g,n}\longrightarrow\mathbb{H}
$ so that $
r\circ\tau=\operatorname{id}_{\mathbb{H}}.$

The cotangent bundle $Q \cT_{g,n}$ is the bundle whose fiber over a Riemann surface $X$ is the collection of all holomorphic quadratic differentials on $X$ where all poles are simple and appear at the punctures. Zeros of the quadratic differentials are also allowed to appear at punctures. The bundle $Q \cT_{g,n}$ admits an action of $\mathrm{GL}(2, \mathbb{R})$ that preserves a stratification given by specifying the orders $\kappa$ of all zeros and poles of the quadratic differential. The corresponding stratum is denoted by $\cQ(\kappa)$. Teichm\"uller showed that Teichm\"uller disks are precisely projections of $\mathrm{GL}(2, \mathbb{R})$-orbits of quadratic differentials to Teichm\"uller space. 

Kra \cite{Kra1981} showed that complex geodesics are holomorphic retracts when they are generated by quadratic differentials all of whose zeros not occurring at punctures have even order. This result was independently discovered by McMullen \cite[Theorem 4.1]{Mc} and we recall its proof in the appendix. Gekhtman--Markovic \cite{GM2020} posited that the converse of Kra's theorem holds and showed that this is true in $\cT_{0,5} \cong \cT_{1,2}$. Their result was extended by Lin--Su \cite{LinSu2026} to $\cT_{0,6} \cong \cT_{2,0}$. We further extend their result to $\cT_{0,n}$ for all $n$.
%Our first main result is the following.

\begin{thm}\label{T:KraGenusZero}
The converse of Kra's Theorem holds in $\cT_{0,n}$.
\end{thm}

Curiously, with one possible exceptional complex geodesic, we find that the same holds in genus one. The \emph{ornithorynque}, introduced by Forni--Matheus \cite{ForniMatheusOrnithorynque}, is a genus one (family of) quadratic differentials in $\cQ(1^3, -1^3)$ formed by pulling back a pillowcase, i.e. a quadratic differential on $\bP^1$ with four poles and no zeros, under a $\bZ/3$ cyclic cover whose branch locus consists of all but one of the poles. 

\begin{thm}\label{T:KraGenusOne}
The converse of Kra's Theorem holds in $\cT_{1,n}$ except possibly for Teichm\"uller disks in $\cQ(1^3, -1^3,0^m)$, for $m \geq 0$, that project to the ornithorynque.
\end{thm}

These results follow from knowing that the converse of Kra's theorem is generically true with five possible exceptions. A component of a stratum is \emph{unmergable} if we cannot merge a simple pole with a zero while remaining in the closure of the component. (See Apisa--Aygun \cite{ApisaAygun} for a definition of merging). 

\begin{thm}\label{T:Generic}
With five possible exceptions, in a component of a stratum $\cQ(\kappa) \subseteq \cQ \cT_{g,n}$ with an odd order zero in $\kappa$ not appearing at a puncture, a full measure set of Teichm\"uller geodesics are not holomorphic retracts. The complement of this full measure set is a locus of covers of genus zero strata. The exceptional components are the unique unmergable component of
$$\cQ(9,-1,0^m), \cQ(6,3,-1,0^m),\cQ(5,-1,0^m), \cQ(4,1,-1,0^m), \text{ or, } \cQ(2,2,1,-1,0^m)$$
where $m$ is a nonnegative integer.
\end{thm}

Finally, we recall that understanding the converse of Kra's theorem is motivated by complex geometry. The \emph{Carath\'eodory metric} on Teichm\"uller space $\cT_{g,n}$ is given by 
\[ d_{Cara}(x,y)=\sup_{f\in\mathrm{Hol}(\cT_{g,n},\mathbb D)}
 d_{hyp}(f(x),f(y)),
\]
where the supremum is taken over all holomorphic maps to $\bD$. Lempert \cite{Lempert1982} showed that on a bounded convex domain in $\bC^n$, the Carath\'eodory and Kobayashi metrics agree. This also holds for unbounded convex domains by \cite[Lemma~3.1]{BracciSaracco2009}. The Kobayashi metric on Teichm\"uller space is the Teichm\"uller metric and the Teichm\"uller disks for which the two metrics coincide are precisely the ones that are holomorphic retracts (see \cite[Lemma~1.3]{GM2020}). Markovic \cite{Markovic2018} used this observation to show that $\cT_{g,n}$ cannot be realized as a convex domain when its dimension is greater than $1$. The noncoincidence of the Carath\'eodory and Kobayashi metrics for Teichm\"uller spaces of infinite type surfaces has been recently established, with seven possible exceptions, by Lin--Markovic \cite[Theorem~1.1]{LinMarkovic2026}.

\subsection*{Organization} In Section \ref{S:Prelims}, we summarize the relevant antecedent work of Kra, McMullen, Markovic, Gekhtman, Apisa, and Wright. In order to be self-contained, all proofs of these results are recalled in the appendix. In Section \ref{S:Proof} we prove Theorem \ref{T:Generic} modulo a technical geometric construction with flat surfaces that we establish in Section \ref{S:Construction}. Deducing Theorems \ref{T:KraGenusZero} and \ref{T:KraGenusOne} from Theorem \ref{T:Generic} comes down to classifying cylinder-free AIS in those genera, which is the focus of Section \ref{S:lowgenus}.

\subsection*{AI Disclosure} All mathematical ideas and writing came from the human authors. ChatGPT was used to create figures and, after the authors had finished revising the paper, to double-check its correctness.

\subsection*{Acknowledgements} Part of this work was carried out at the Institute for Advanced Study, Princeton, during the Workshop on Several Complex Variables and Related Topics. We are grateful to the Institute for its hospitality and to the Minerva Research Foundation for its generous support of the workshop. The arguments in Section \ref{S:Proof} were developed in discussions with and with input from Bernhard Lamel and Rumpa Masanta.

P.A. was partially supported by a Sloan Research Fellowship and NSF grant no.~DMS-2304840. J.A. was partially supported by the Simons Dissertation
Fellowship. D.N. was supported by the Kovner Member Fund, by the Israel Science Foundation grant 1167/17 and by Minerva grant 714141. A.Z. was partially supported by a Sloan Research Fellowship and grants DMS-2105580 and DMS-2452068 from the National Science Foundation.

\section{Preliminaries} \label{S:Prelims}

Given a component $\cC$ of a stratum of quadratic differentials, let $\cM$ be the locus quadratic differentials inside $\cC$ which generate Teichm\"uller disks, also called \emph{complex geodesics}, that are holomorphic retracts of Teichm\"uller space.

\begin{prop}\cite[Proposition 2.1]{GM2020}
$\cM$ is an affine invariant subvariety.
\end{prop}
\begin{proof}
By Eskin, Mirzakhani, Mohammadi \cite{EM, EMM}, it suffices to show that $\cM$ is closed and $\mathrm{GL}(2, \bR)$ invariant. Invariance is immediate so suppose that $x_n$ is a sequence of quadratic differentials in $\cM$ converging to a point $x$ in $\cC$. We must show that $x \in \cM$, i.e. that its Teichm\"uller disk is a holomorphic retract of Teichm\"uller space. Let $f_n$ be retracts of Teichm\"uller space to the Teichm\"uller disk generated by $x_n$ which sends $x_n$ to $i$ after identifying the disk with $\mathbb{H}$. Then $f_n$ is a sequence of contractions and so, by Arzel\`a-Ascoli, a subsequence converges to the desired retraction to the disk generated by $x$.
\end{proof}

A quadratic differential $q$ on a surface $X$ is \emph{horizontally periodic} if its nonsingular horizontal trajectories are compact. In this case, $(X,q)$, which can be thought of as a half-translation surface, is a union of horizontal cylinders $C_1,\dots, C_k$. For $(\lambda_1,\dots, \lambda_k) \in \bH^k$, let  $\tilde \cE_X(\lambda_1,\dots, \lambda_k)$ be the quadratic differential obtained by applying the affine map 
$$
x+yi  \mapsto x+\lambda_j y
$$
to the $j^{th}$ cylinder. The image of this map is called the \emph{polyplane}. Let $\cE_X(\lambda_1,\dots, \lambda_k)$ be the image of $\tilde \cE_X(\lambda_1,\dots, \lambda_k)$ in $\cT_{g,n}$. So 
$$
\lambda \in \bH \mapsto \cE_X(\lambda, \dots, \lambda)
$$
parametrizes the Teichm\"uller disk generated by $(X,q)$. 
The following two result were originally shown by Markovic \cite{Markovic2018} for $\cT_{0,5}$ and extended to $\cT_{g,n}$ by Gekhtman--Markovic \cite{GM2020}. 

\begin{thm}\cite[Theorem 4.1]{GM2020}\label{thm:retraction on polyplane} Let $X$ be a unit-area horizontally periodic quadratic differential with horizontal cylinders of area $\alpha_1, \hdots, \alpha_k$. Suppose that its Teichm\"uller disk is a holomorphic retract of Teichm\"uller space. Then the retract $\Phi$ may be chosen so that $\Phi\cE_X(\lambda_1,\dots, \lambda_k) = \sum \alpha_j \lambda_j$.
\end{thm} 

The proof of Theorem \ref{thm:retraction on polyplane} is recalled in Appendix \ref{A:Background}.

\begin{thm}\cite[Corollary 4.4]{GM2020}\label{thm:polyplanes are all or nothing} If $X$ is horizontally periodic and the associated complex geodesic admits a holomorphic retract, then so does $\tilde \cE_X(\lambda_1,\dots, \lambda_k)$ for any  $(\lambda_1,\dots, \lambda_k) \in \bH^k$.
\end{thm}
\begin{proof}
Another point $x'$ in the polyplane can be described as applying the matrices $\begin{pmatrix} 1 & c_i \\ 0 & d_i \end{pmatrix}$ to the $i$th cylinder for some choice of $c_i$ and $d_i$ for each $i$. Let $\Phi$ and $(\alpha_i)_{i=1}^k$ be as in Theorem \ref{thm:retraction on polyplane}. Let $d = \sum_i \alpha_i d_i$ and $c = \sum_i \alpha_i c_i$. Then $\frac{\Phi - c}{d}$ is the desired retraction to the Teichm\"uller disk generated by $x'$.
\end{proof}

%After the transformation, areas are d_i a_i. We want (d_i lambda_i + c_i) to be sent to i. Under Phi it goes to sum_i a_i(d_i lambda_i + c_i). So (Phi - c)/d goes to sum_i a_i lambda_i 
%
% Transformations send lambda to (x_i+c)+(d_iy). Phi sends this to sum_i a_id
%

An invariant subvariety $\cM$ is called \emph{cylinder-free} if every cylinder on a surface in $\cM$ can be arbitrarily sheared and dilated while fixing its complement and while remaining in $\cM$.

\begin{cor}\label{C:CylinderFree}
The locus of retracts $\cM$ is cylinder-free.
\end{cor}
\begin{proof}
By Minsky--Weiss \cite{MinW}, $\cM$ contains a horizontally periodic quadratic  differential. By Wright \cite{Wcyl}, if the polyplane associated to every horizontally periodic differential in $\cM$ is contained in $\cM$, then $\cM$ is cylinder-free. This latter condition is guaranteed by Theorem \ref{thm:polyplanes are all or nothing}.
\end{proof}

An interesting problem is to classify the cylinder-free invariant subvarieties. This is essentially a topological problem by the following result.

\begin{thm}\cite[Theorem 1.1]{ApisaWrightGeminal}\label{T:AW}
The only nonempty proper cylinder-free invariant subvarieties are full loci of covers of strata of genus zero quadratic differentials.
\end{thm}

%By Theorem \ref{T:AW}, to prove Theorem \ref{T:Generic}, it suffices to find just one quadratic differential in a stratum with an odd order zero which does not generate a disk which admits a holomorphic retract.

%%%%%%%%
%
% New section
%
%%%%%%%%%

\section{Proof of Theorem \ref{T:Generic} given Theorem \ref{T:JSsurface}} \label{S:Proof}

In the sequel, we will consider a component $\cC$ of a stratum $\cQ \subseteq T^*\cT_{g,n}$ of quadratic differentials. All of the poles of the quadratic differentials must be marked and all other marked points must occur at zeros (possibly at zeros of order zero). However, not all zeros are necessarily marked. Nonetheless, it will be useful to consider the closure $\overline{\cQ}$ of $\cQ$ in the multiscale compactification \cite{BCGGM} where every pole and zero is marked.

%\subsection{Proof of Theorem \ref{T:Generic} given Theorem \ref{T:JSsurface}}

\begin{figure}
\begin{tikzpicture}[
  x=1cm,y=1cm,
  line width=.65pt,
  line cap=round,line join=round,
  font=\small,
  surface/.style={draw=black,fill=blue!5},
  core/.style={densely dashed,line width=.5pt},
  brace/.style={
    decorate,
    decoration={brace,amplitude=4pt},
    line width=.5pt
  }
]

% ============================================================
% Upper left: Q(1,-1^5), flat model.
% ============================================================
\begin{scope}
  \node at (1.35,2.55) {$\mathcal Q(1,-1^5)$};

  \path[surface]
    (0,0)--(2.7,0)--(2.7,1)
    --(1.05,1)--(1.05,2)--(0,2)--cycle;

  \draw[core] (0,1.52)--(1.05,1.52);
  \draw[core] (0,.43)--(2.7,.43);
  \node[below] at (.57,1.52) {$q_0$};
  \node[below] at (1.73,.43) {$q_1$};

  \draw[brace] (-.19,1)--(-.19,2)
    node[midway,left=6pt] {$a_0$};
  \draw[brace] (-.19,0)--(-.19,1)
    node[midway,left=6pt] {$a_1$};

  \ticks{0}{1.65}{90}{1}
  \ticks{1.05}{1.65}{90}{1}
  \ticks{0}{.67}{90}{2}
  \ticks{2.7}{.67}{90}{2}

  \chevron{.30}{2}{0}
  \chevron{.77}{2}{180}
  \ticks{1.43}{1}{0}{1}
  \chevron{1.55}{1}{0}
  \ticks{2.27}{1}{0}{1}
  \chevron{2.15}{1}{180}
  \ticks{.40}{0}{0}{2}
  \chevron{.55}{0}{0}
  \ticks{2.15}{0}{0}{2}
  \chevron{1.99}{0}{180}

  \crosspt{0}{2}
  \crosspt{.54}{2}
  \crosspt{1.86}{1}
  \crosspt{1.35}{0}
  \crosspt{0}{0}
  \dotpt{0}{1}
  \dotpt{1.05}{1}
  \dotpt{2.7}{1}
\end{scope}

% Upper left: corresponding sphere.
\begin{scope}[shift={(4.45,1.02)}]
  \sphere{0}{0}{.87}{.94}
  \foreach \x/\y in {
    -.61/.35,-.32/.53,0/.60,.32/.53,.61/.35
  }{\crosspt{\x}{\y}}
  \dotpt{0}{-.94}
  \node[below=3pt] at (0,-.94) {$1$};
\end{scope}

% ============================================================
% Upper right: H(2,0^5), flat model.
% ============================================================
\begin{scope}[shift={(6.25,0)}]
  \node at (1.25,2.55) {$\mathcal H(2,0^5)$};

  \path[surface]
    (0,0)--(2.45,0)--(2.45,1)
    --(1.05,1)--(1.05,2)--(0,2)--cycle;

  \draw[core] (0,1.38)--(1.05,1.38);
  \draw[core] (0,.30)--(2.45,.30);
  \node[below] at (.58,1.38) {$\widetilde q_0$};
  \node[above] at (1.31,.30) {$\widetilde q_1$};

  \draw[brace] (-.19,1)--(-.19,2)
    node[midway,left=6pt] {$\widetilde a_0$};

  \ticks{1.05}{1.53}{90}{1}
  \ticks{0}{.57}{90}{2}
  \ticks{2.45}{.57}{90}{2}
  \ticks{.53}{2}{0}{3}
  \ticks{.53}{0}{0}{3}
  \ticks{2.06}{1}{0}{4}
  \ticks{2.06}{0}{0}{4}

  \foreach \x/\y in {
    0/2,1.05/2,0/1,1.05/1,2.45/1,
    0/0,1.12/0,2.45/0
  }{\dotpt{\x}{\y}}
  \foreach \x/\y in {
    0/1.65,.53/1.65,1.76/1,.55/.68,1.84/.68
  }{\crosspt{\x}{\y}}
\end{scope}

% Upper right: smooth genus-two surface.
\begin{scope}[shift={(10.45,1.05)}]
  \path[surface]
    (0,-.95)
    .. controls (-.84,-.94) and (-.88,-.22) .. (-.64,.42)
    .. controls (-.45,1.03) and (.41,1.17) .. (.65,.55)
    .. controls (.96,-.20) and (.80,-.84) .. (0,-.95);
  \handle{0}{.30}
  \handle{0}{-.20}
  \crosspt{.13}{.96}
  \crosspt{-.07}{.43}
  \crosspt{.06}{.23}
  \crosspt{-.07}{-.07}
  \crosspt{.06}{-.27}
  \dotpt{0}{-.95}
  \node[below=3pt] at (0,-.95) {$2$};
\end{scope}

% ============================================================
% Lower left: Q(0,-1^4), flat model.
% ============================================================
\begin{scope}[shift={(0,-2.65)}]
  \node at (1.775,1.65) {$\mathcal Q(0,-1^4)$};

  \path[surface] (0,0) rectangle (2.7,1);
  \draw[core] (0,.36)--(2.7,.36);
  \node[below] at (1.83,.46) {$q_1$};

  \draw[brace] (0,1.20)--(1.35,1.20)
    node[midway,above=6pt] {$q_0$};
  \draw[brace] (-.19,0)--(-.19,1)
    node[midway,left=6pt] {$a_1$};

  \ticks{0}{.65}{90}{2}
  \ticks{2.7}{.65}{90}{2}
  \chevron{.39}{1}{0}
  \chevron{1.05}{1}{180}
  \ticks{1.62}{1}{0}{1}
  \chevron{1.73}{1}{0}
  \ticks{2.34}{1}{0}{1}
  \chevron{2.23}{1}{180}
  \ticks{.50}{0}{0}{2}
  \chevron{.65}{0}{0}
  \ticks{2.10}{0}{0}{2}
  \chevron{1.94}{0}{180}

  \openpt{0}{1}
  \openpt{1.35}{1}
  \openpt{2.7}{1}
  \crosspt{.72}{1}
  \crosspt{1.98}{1}
  \crosspt{1.35}{0}
\end{scope}

% Lower left: two spherical components.
\begin{scope}[shift={(4.45,-1.81)}]
  \sphere{0}{0}{.67}{.57}
  \foreach \x/\y in {
    -.43/.22,-.16/.35,.16/.35,.43/.22
  }{\crosspt{\x}{\y}}

  \sphere{0}{-.94}{.35}{.37}
  \dotpt{0}{-.57}
  \dotpt{-.23}{-1.22}
  \crosspt{.23}{-1.22}
  \node[below left=2pt] at (-.23,-1.22) {$1$};
  \node[below right=2pt] at (.23,-1.22) {$-1$};
\end{scope}

% ============================================================
% Lower right: H(0^6), flat model.
% ============================================================
\begin{scope}[shift={(6.25,-2.65)}]
  \node at (1.75,1.65) {$\mathcal H(0^6)$};

  \path[surface] (0,0) rectangle (2.45,1);
  \draw[core] (0,.32)--(2.45,.32);
  \node[above] at (1.27,.22) {$\widetilde q_1$};

  \draw[brace] (0,1.20)--(1.20,1.20)
    node[midway,above=6pt] {$\widetilde q_0$};

  \ticks{0}{.55}{90}{2}
  \ticks{2.45}{.55}{90}{2}
  \ticks{.56}{1}{0}{3}
  \ticks{.56}{0}{0}{3}
  \ticks{1.91}{1}{0}{4}
  \ticks{1.91}{0}{0}{4}

  \openpt{0}{1}
  \openpt{0}{0}
  \squarept{1.20}{1}
  \squarept{1.20}{0}
  \crosspt{.72}{1}
  \crosspt{2.07}{1}
  \crosspt{.63}{.61}
  \crosspt{1.79}{.61}
\end{scope}

% Lower right: nodal surface.
\begin{scope}[shift={(10.45,-1.75)}]
  % Lower component.
  \path[surface]
    (-.40,-.37)
    .. controls (-.63,-.69) and (-.63,-1.45) .. (0,-1.50)
    .. controls (.63,-1.45) and (.63,-.69) .. (.40,-.37)
    .. controls (.21,-.47) and (-.21,-.47) .. (-.40,-.37);

  % Upper genus-one component.
  \path[surface]
    (-.40,-.37)
    .. controls (-.83,-.10) and (-.55,.63) .. (0,.65)
    .. controls (.55,.63) and (.83,-.10) .. (.40,-.37)
    .. controls (.17,-.52) and (-.17,-.52) .. (-.40,-.37);

% Unshaded space between the torus and the rational component.
\path[fill=white,draw=none]
  (-.40,-.37)
  .. controls (-.17,-.52) and (.17,-.52) .. (.40,-.37)
  .. controls (.41,-.91) and (-.41,-.91) .. (-.40,-.37)
  -- cycle;

% Restore the lower boundary of the torus.
\draw (-.40,-.37)
  .. controls (-.17,-.52) and (.17,-.52) .. (.40,-.37);

  % The two branches joining the components.
  \draw (-.40,-.37)
    .. controls (-.41,-.91) and (.41,-.91) .. (.40,-.37);

  \handle{0}{.15}
  \crosspt{.08}{.65}
  \crosspt{-.07}{.28}
  \crosspt{.06}{.08}
  \crosspt{0}{-.48}
  \crosspt{0}{-.77}

  \openpt{-.40}{-.37}
  \openpt{.40}{-.37}
  \dotpt{0}{-1.50}
  \node[below=3pt] at (0,-1.50) {$2$};
\end{scope}

\end{tikzpicture}
    \caption{The periods of $a_0, a_1, q_0, q_1$ provide coordinates on $\cQ(1,-1^5) \subseteq T^* \cT_{0,5}$. Degenerating the top cylinder causes $a_0$ to vanish and is drawn in a flat and multiscale perspective as the bottom left surface(s). If the simple zero is not marked, the degeneration does not pass to the boundary of $\cT_{0,5}$. But the same is not true for the corresponding degeneration of the holonomy cover, shown in the right column.}
    \label{F:Coordinates}
\end{figure}
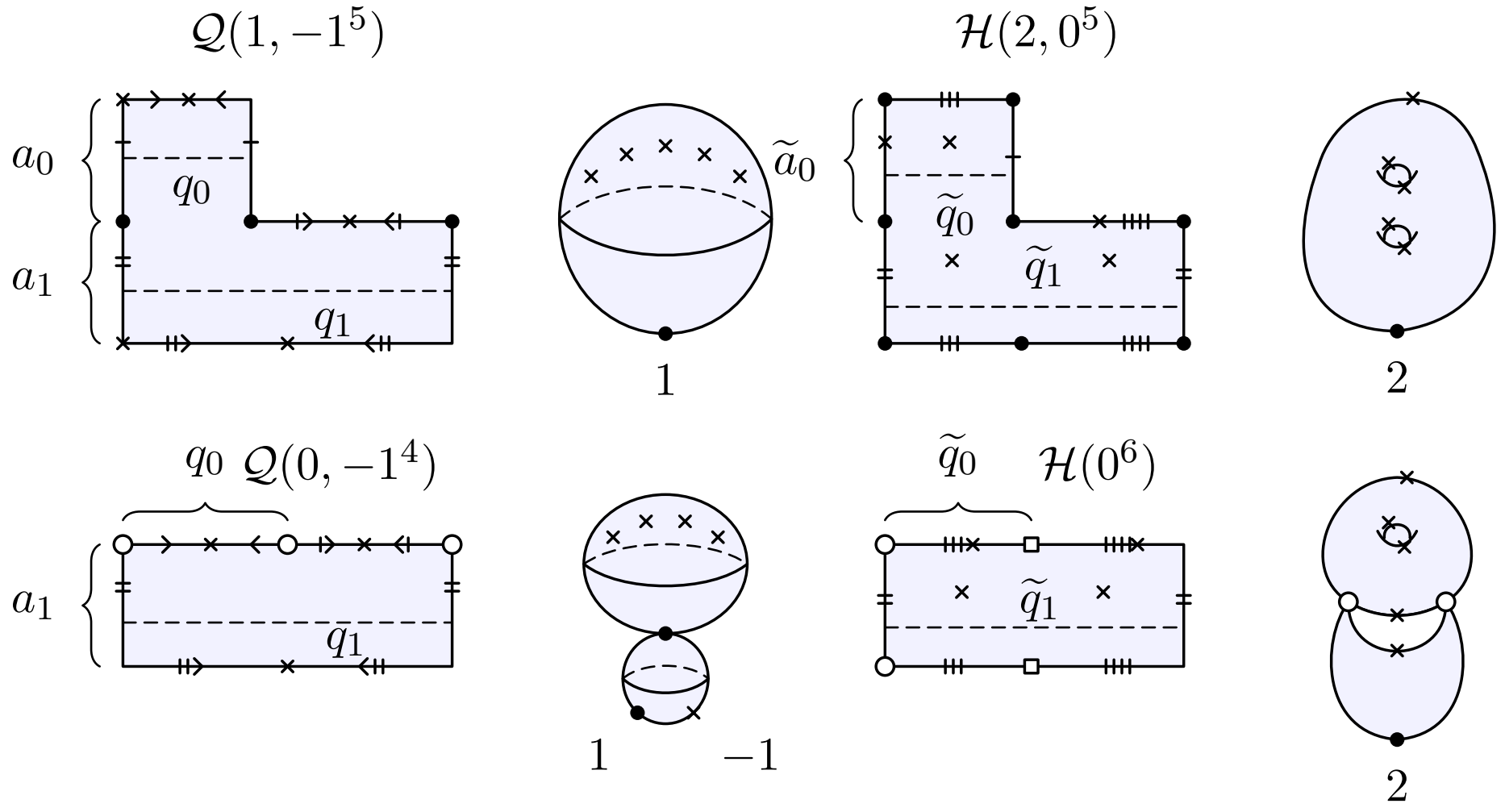

Two singularities $s_1$ and $s_2$ can be \emph{merged} in a component $\cC$ of $\cQ$ if the multiscale compactification of $\cC$ contains a point with exactly two irreducible components, one a thrice-punctured sphere with punctures corresponding to $s_1$ and $s_2$ (see the introduction of \cite{ApisaAygun} for discussion. See Figure \ref{F:Coordinates}.). For example, if $s_1$ and $s_2$ are on opposite boundaries of a cylinder and are joined by a cross curve $\gamma$ which is the only vanishing cycle when the cylinder is collapsed, then $s_1$ and $s_2$ can be merged. (See the left column of Figure \ref{F:Coordinates} for an illustration and \cite{ApisaWrightDiamonds} for a discussion of cylinder collapses).

 %For the rest of this paper, we will abuse notation and call a singularity we refer to by its order.
In the sequel, we will often let $z$ refer to both a zero and its order. The proof of Theorem \ref{T:Generic} will rely on the following theorem, whose proof we defer. 

%Two cylinder core curves are \emph{hat homologous} if they remain parallel to each other on all nearby surfaces in the stratum.

\begin{thm} \label{T:JSsurface} 
If $z$ is an odd order zero and $\cC$ is not an unmergeable component of $$\cQ(9,-1,0^m), \cQ(6,3,-1,0^m),\cQ(5,-1,0^m), \cQ(4,1,-1,0^m), \cQ(2,2,1,-1,0^m)$$ for some $m\in\mathbb{Z}_{\geq 0}$, then $\cC$ contains a horizontally periodic surface  with horizontal cylinders $(C_i)_{i=0}^{k-1}$ with core curves $(q_i)_{i=0}^{k-1}$ and cross curves $(a_j)_{j=0}^{k-1}$ with the following properties: 
\begin{enumerate}
    \item Collapsing $C_0$ sends the length of $a_0$ to zero and merges $z$ and a distinct odd order singularity, 
    \item $k$ is either the maximal number $M$ of horizontal cylinders that a periodic surface in $\cQ$ can have or $k = M-1$,
    \item When $k=M$, $(a_i)_{i=0}^{k-1} \cup (q_i)_{i=0}^{r-1}$ is a basis of period coordinates where $r$ is the rank of $\cQ$,
    \item When $k=M-1$, there is a saddle connection $q_{-1}$ that does not intersect $a_0$ (even at endpoints) and so $(a_i)_{i=0}^{k-1} \cup (q_i)_{i=-1}^{r-1}$ is a basis of period coordinates
    \item The polyplane near the surface is cut out by fixing the periods of each $q_i$. 
    %and the equations $\emph{Im}(a_i) >0$. {\color{red} added last part}
    %The surface contains a horizontal cylinder $C$, with core curve $q_0$ and cross curve $a_0$, which, when collapsed, degenerates $a_0$ and merges $z$ with a distinct odd order singularity, and    \item There is a collection of curves $q_0\cup (a_j)_{j=0}^{h-1} \cup (q_i)_{i=1}^{2g+n-h-3}$ whose periods give local coordinates for $\cQ$ such that the polyplane is cut out by fixing the periods of $q_0$ and each $q_i\neq q_0$ pairs trivially with the vanishing cycles of the collapse of $C$.
\end{enumerate} 
\end{thm}
Note that the strata $\cQ(5,-1), \cQ(4,1,-1), \cQ(2,2,1,-1)$ each consists of only one component, and this component is unmergeable.

Let $(X, q)$ be a unit area surface produced by Theorem \ref{T:JSsurface} and take $z$ to be an unmarked zero. Let $(Y, \omega)$ be its holonomy cover. As in the theorem, $a_i$ and $q_i$ will be used for both curves and their periods, which are taken on the holonomy cover and hence well-defined up to sign. Let $\mathrm{Col}_{C_0}(X,q)$ be the result of collapsing the cylinder $C_0$. 

%Let $\wt{C}$ be the lift of $C$ to the holonomy cover $(Y, \omega)$ of $(X, q)$. 

Smooth points in $\cQ$ near $\mathrm{Col}_{C_0}(X,q)$ are formed by plumbing the node, i.e. replacing the node, which is modeled on solutions to $\{ (x,y) \in \bC^2 : xy = 0\}$ by the annulus $\{(x, y) \in \bC^2 : xy = t \}$. The parameter $t$ is called a \emph{plumbing parameter}, which is holomorphic on a neighborhood of $\mathrm{Col}_{C_0}(X, q)$ in $\overline{\cQ}$.

%(and thus the cross curve $a_0$), produced by Theorem \ref{T:JSsurface} which merges $z$ and some odd order singularity (the operation of cylinder collapses is discussed in detail in \cite{ApisaWrightDiamonds}). Using log periods, defined below, which essentially encode the Picard-Lefschetz formula, we can extend $q_0\cup (a_j)_{j=0}^{h-1} \cup (q_i)_{i=1}^{2g+n-h-3}$ to a coordinate system of $\overline{\cQ}$ at $(X', q')$. 

Collapsing $C_0$ pinches a separating curve $\gamma_v$ that pinches off a copy of $\bP^1$ containing $z$ and the singularity $w$ with which it merges. On the holonomy cover, $\gamma_v$ has two lifts and pinching them pinches off a copy of $\bP^1$ that contains the preimages of $z$ and $w$. See Figure \ref{F:Coordinates}. Call the two (positively-oriented) vanishing cycles on $(Y, \omega)$, $\eta_1$ and $\eta_2$.  Suppose that $\gamma$ is the period of any relative cycle on $(Y, \omega)$ that is not contained in the nodal $\bP^1$. Benirschke \cite{MR4662298} defines its \emph{log period} at points near $\mathrm{Col}_{C_0}(X, q)$ as
\[ \psi_\gamma := \int_\gamma \omega + \langle \gamma, \eta_1 - \eta_2 \rangle \log t\int_{\eta_1} \omega\]
where $\langle \cdot, \cdot \rangle$ is the intersection pairing and $t$ is the plumbing parameter. By Benirschke \cite[Theorem 5.2]{MR4662298}, log periods are holomorphic in a neighborhood of $\mathrm{Col}_{C_0}(X, q)$ in $\overline{\cQ}$.  Letting a prime denote a log period, for each $i$ and some $\epsilon_i \in \{0,\pm1, \pm 2\}$,
% \[ a_i' := a_i + \epsilon_i a_0 \log t = a_i + \frac{\epsilon_i}{1+\kappa} a_0 (\log a_0 + f) \]
\[ a_i' := a_i + \epsilon_i a_0 \log t. \]

%In our setting, the degeneration of $(X, q)$ is always a merge and the degeneration of its holonomy cover involves pinching two homologous non-separating simple closed curves that form the boundary of a $\bP^1$ that is punctured at two points, i.e. the preimages of the two odd order singularities that merge on $(X, q)$. Call the two (positively-oriented) vanishing cycles $\gamma_1, \gamma_2$. Points near $(X', q')$ are formed by plumbing the node, i.e. replacing the node, which is modeled on solutions to $\{ (x,y) \in \bC^2 : xy = 0\}$ by the annulus $\{(x, y) \in \bC^2 : xy = t \}$. The parameter $t$ is called a \emph{plumbing parameter}, which is holomorphic near the boundary of $\overline{\cQ}$. Suppose that $\gamma$ is the period of any relative cycle that is not contained in the nodal $\bP^1$. Benirschke \ann{J: Should we get rid of the discussion of this in the previous section.} defines its \emph{log period} at points near $(X', q')$ as

%where $f$ is some nonvanishing holomorphic function. 

\begin{lem}\label{L:HolomorphicCoordinates}
The plumbing parameter and the log periods of the $q_i$'s and $(a_i)_{i=1}^{k-1}$ are a holomorphic coordinate system for $\overline{\cQ}$ near $\mathrm{Col}_{C_0}(X, q)$. The same holds if each $a_i'$ is replaced by $a_i'' := a_i - \frac{\epsilon_i}{2} q_0$.
%
%The functions $(q_0', a_1', \hdots, a_{h-1}',q_1', \hdots, q'_{2g+n-h-3},t)$ provide a holomorphic coordinate system for $\overline{\cQ}$ near $(X', q')$. Likewise is true if each $a_i'$ is replaced by $a_i'' := a_i - \frac{\epsilon_i}{2} q_0$.
\end{lem}

%{\color{red} A: Is it worth saying something in the proof below like, BCGGM construct holomorphic coordinates on the bundle over Moduli space but since we are in the interior it lifts to the bundle over Teichm\"uller space?}

\begin{proof}
The coordinates (less $t$) provide period coordinates on the boundary stratum of $\overline{\cQ}$ containing $\mathrm{Col}_{C_0}(X, q)$. To see this, note that when $t = 0$, the log periods are the usual periods in the WYSIWYG compactification and these periods form a basis. By \cite[Theorems 10.1, 10.2]{BCGGM}, adding the plumbing coordinate $t$ creates a holomorphic coordinate system\footnote{Their result is established over $\overline{\cM}_{g,n}$ rather than $\overline{\cT}_{g,n}$. However, since $\mathrm{Col}_{C_0}(X,q)$ lies over the interior of $\overline{\cT}_{g,n}$, and hence $\overline{\cM}_{g,n}$ by forgetting the marking, we may choose a marking of the underlying family in a sufficiently small neighborhood and obtain a holomorphic local lift of these coordinates to $\overline{\cT}_{g,n}$.}. The same argument applies to the second proposed coordinate system, i.e.\ it forms a period coordinate system on the boundary stratum and hence adding the plumbing parameter extends it to a coordinate system on $\overline{\cQ}$ by the inverse function theorem. Since $a_i'' = a_i' - \frac{\epsilon_i}{2} q_0'$, it is also holomorphic. 
\end{proof}

Say the \emph{extended polyplane} is the locus containing $(X, q)$ where the periods of all $q_i$'s are fixed except for that of $q_0$. One key point is that since all $q_i$'s except $q_0$ avoid $a_0$, the log period of such a $q_i$ is just its period. So, the extended polyplane is analytic even at the boundary of its parametrization by $\tilde \cE_X$. The fact that the solution set of $q_0 = c$ in the extended polyplane for some constant $c$ is not analytic will be key. %Whereas, the extended polyplane is analytic because each lift of $q_i \neq q_0$ pairs trivially with the vanishing cycles on the holonomy cover, and hence $q_i' = q_i$ (i.e., setting $q_i$ equal to a constant defines an analytic condition even at the boundary). Another key point is that fixing $(q_i)_{i=1}^{2g+n-h-3}$ and varying $q_0$ only by real numbers defines horizontally periodic surfaces in a neighborhood of $(X, q)$. 

\begin{cor}\label{C:MarkovicCriterion}
Suppose that the Teichm\"uller disk generated by $(X, q)$ is a holomorphic retract of Teichm\"uller space and that the period of $a_j$ is $\sqrt{-1}$ for all $j$. Let $c$ be the period of $q_0$ and $\alpha_i$ the area of cylinder with cross curve $a_i$ on $(X, q)$. Then the retract may be taken to be $b(c-q_0) + \sum_i \alpha_i a_i$ on the extended polyplane where $b = \sum_i \frac{\alpha_i \epsilon_i}{2}$.
% Then on the extended polyplane there are real constants $\alpha_i, \beta_i$ so that the retract is given by $\beta_{-1} + \sum_i \beta_i a_i''$. Near $(X,q)$ this is given by $\sum_i \alpha_i a_i$.
\end{cor}
\begin{proof}
By Markovic's criterion (Theorem \ref{thm:retraction on polyplane}), the holomorphic retract on the polyplane may be given by $\sum_i \alpha_i a_i = bq_0 + \sum_i \alpha_i a_i''$ since each $a_i = \sqrt{-1}$ on $(X,q)$. The polyplane containing $(X, q)$ is defined within the extended polyplane by $q_0 = c$. Let $H$ be the holomorphic retract restricted to the extended polyplane. Then $H - (bc) - \sum_i \alpha_i a_i''$ is holomorphic and vanishes on the locus $q_0 = c$. In the holomorphic coordinate system given by Lemma \ref{L:HolomorphicCoordinates}, this is the locus defined by $q_0' = c + 2 f t^{\mu+1} \log t $ where $\mu$ is the order (at the top level) of a nodal singularity on the holonomy cover and $f$ is the nonvanishing holomorphic function given by $a_0/t^{\mu+1}$ (see \cite[Section 11]{BCGGM} for a proof of the these properties of $f$). This locus is nonanalytic as can be seen by replacing the coordinate $q_0'$ with $q_0'' := \frac{q_0-c}{2f}$ and considering the locus given by $q_0'' = t^{\mu+1} \log t$. So $q_0 = c$ is a nonanalytic codimension one subset of the extended polyplane on which the (analytic) equation $H - (bc) - \sum_i \alpha_i a_i'' = 0$ holds. So the equation holds on the entire extended polyplane, i.e. $H = (bc) + \sum_i \alpha_i a_i''$. 
\end{proof}

The following is a restatement and proof of Theorem \ref{T:Generic}.

\begin{thm}\label{T:GenericKra}
The converse of Kra's Theorem holds generically on all strata except those excluded in the hypotheses of Theorem \ref{T:JSsurface}.
\end{thm}

%{\color{red} A: I find it kind of annoying that in the previous proof the $q_i$ and $a_j$ are functions, but not a big deal... we could say something like Let $(q_{0,i})_{i=0, -1}^{r-1}$ and $(a_{0,j})^{k-1}_{j = 0}$ be the values of $(q_i)_{i=0, -1}^{r-1}$ and $(a_j)^{k-1}_{j = 0}$ at $(X,q)$.

%Along similar lines, it is kind of unpleasant that the $q_i$ are periods but $q$ is the quadratic differential. We could do something like use the mathscr package and use mathscr{q} for the quadratic differential. But again not really a big deal.

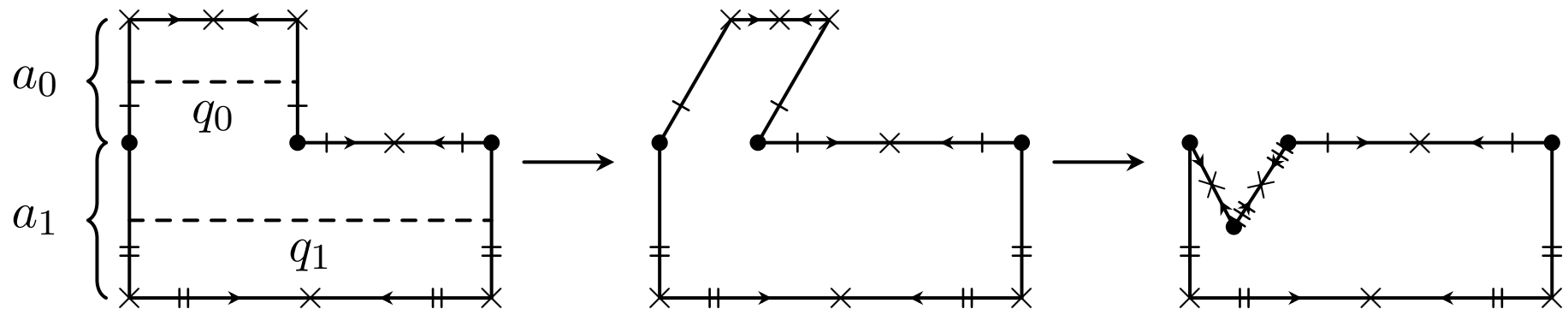
\begin{figure}[H]
\centering

% Figure:
\begin{tikzpicture}[
  line width=.75pt,
  line cap=round,
  line join=round,
  font=\small
]

% ============================================================
% 1. Original surface
% ============================================================
\begin{scope}
  \coordinate (A) at (0,0);
  \coordinate (B) at (2.8,0);
  \coordinate (C) at (2.8,1.2);
  \coordinate (D) at (1.3,1.2);
  \coordinate (E) at (1.3,2.15);
  \coordinate (F) at (0,2.15);
  \coordinate (G) at (0,1.2);

  % Cylinder sides: ticks only, no arrowheads.
  \sideedge{A}{G}{2}
  \sideedge{B}{C}{2}
  \sideedge{G}{F}{1}
  \sideedge{D}{E}{1}

  \foldedge{F}{E}{P}{0}
  \foldedge{D}{C}{Q}{1}
  \foldedge{A}{B}{R}{2}

  \draw[dashed]
    (0,1.67)--node[below] {$q_0$}(1.3,1.67);
  \draw[dashed]
    (0,.60)--node[below] {$q_1$}(2.8,.60);

  \draw[decorate,decoration={brace,amplitude=4pt}]
    (-.18,1.2)
    --node[left=6pt] {$a_0$}(-.18,2.15);
  \draw[decorate,decoration={brace,amplitude=4pt}]
    (-.18,0)
    --node[left=6pt] {$a_1$}(-.18,1.2);

  \foreach \p in {A,B,E,F}{\pole{\p}}
  \edgepole{P}{F}{E}
  \edgepole{Q}{D}{C}
  \edgepole{R}{A}{B}

  \foreach \p in {C,D,G}{\zeropoint{\p}}
\end{scope}

\draw[-{Stealth[length=4pt,width=4pt]}]
  (3.05,1.05)--(3.75,1.05);

% ============================================================
% 2. Small circumference
% ============================================================
\begin{scope}[xshift=4.1cm]
  \coordinate (A) at (0,0);
  \coordinate (B) at (2.8,0);
  \coordinate (C) at (2.8,1.2);
  \coordinate (D) at (.76,1.2);
\coordinate (E) at (1.31,2.15);
\coordinate (F) at (.55,2.15);
  \coordinate (G) at (0,1.2);

  \sideedge{A}{G}{2}
  \sideedge{B}{C}{2}
  \sideedge{G}{F}{1}
  \sideedge{D}{E}{1}

  \foldedge{F}{E}{P}{0}
  \foldedge{D}{C}{Q}{1}
  \foldedge{A}{B}{R}{2}

  \foreach \p in {A,B,E,F}{\pole{\p}}
  \edgepole{P}{F}{E}
  \edgepole{Q}{D}{C}
  \edgepole{R}{A}{B}

  \foreach \p in {C,D,G}{\zeropoint{\p}}
\end{scope}

\draw[-{Stealth[length=4pt,width=4pt]}]
  (7.15,1.05)--(7.85,1.05);

% ============================================================
% 3. Overcollapsed surface
% ============================================================
\begin{scope}[xshift=8.2cm]
  \coordinate (A) at (0,0);
  \coordinate (B) at (2.8,0);
  \coordinate (C) at (2.8,1.2);
\coordinate (D) at (.76,1.2);
\coordinate (E) at (.34,.55);
  \coordinate (F) at (0,1.2);

  \sideedge{A}{F}{2}
  \sideedge{B}{C}{2}

  \foldedge{F}{E}{P}{0}
  % Three ticks distinguish this new fold from the other folds.
  \foldedge{E}{D}{S}{3}
  \foldedge{D}{C}{Q}{1}
  \foldedge{A}{B}{R}{2}

  \foreach \p in {A,B}{\pole{\p}}

  % Rotate the crosses relative to their respective edges.
  \edgepole{P}{F}{E}
  \edgepole{S}{E}{D}
  \edgepole{Q}{D}{C}
  \edgepole{R}{A}{B}

  \foreach \p in {C,D,E,F}{\zeropoint{\p}}
\end{scope}

\end{tikzpicture}

\caption{In the proof of Theorem \ref{T:GenericKra}, we decrease the circumference $q_0$ of $C_0$, slightly shear $C_0$ by $x$, and then overcollapse it.}
\label{F:Overcollapse}
\end{figure}

%} 
\begin{proof}
Suppose, in order to derive a contradiction, that every surface in the component of the stratum with an odd order zero not occurring at a puncture (and to which Theorem \ref{T:JSsurface} applies) generates a Teichm\"uller disk that is a holomorphic retract. Start with $(X, q)$ from Theorem \ref{T:JSsurface}. We may assume that the period of $a_j$ is $\sqrt{-1}$ for all $j$ and that the periods of $(q_i)_{i=0, -1}^{r-1}$ are independent over $\bQ$. Let $F$ be the retract from Corollary \ref{C:MarkovicCriterion}. Now, move to a different surface in the extended polyplane of $(X,q)$ to make the period of $q_0$ sufficiently small and ensure that the periods of $(a_j)^{k-1}_{j = 0}$ are purely imaginary and independent over $\bQ$. We do this while continuing to insist that the periods of $(q_i)_{i=0, -1}^{r-1}$ are independent over $\bQ$. Then shear $C_0$ so that the period of $a_0$ becomes $a_0+x$ where $x$ is some real number that is linearly independent from all periods of $a_i$'s and $q_i$'s. Starting at this surface, consider the path $(X_t)_{t \geq 0}$ where all periods of $(q_i)_{i=0, -1}^{r-1}$ and $(a_j)^{k-1}_{j = 0}$ are fixed except for $a_0(t) = a_0 + x-\sqrt{-1}t$. (See Figure \ref{F:Overcollapse}.) Along this path, the area at time $t$ is $A(t) := q_0 (\text{Im}(a_0) - t) + \sum_{i \ne 0} q_i \text{Im}(a_i)$. By making $q_0$ small beforehand, we can travel along this path for arbitrarily long\footnote{We can travel in the stratum as long as area is bounded away from zero and no period vanishes; since periods are linear combinations of the $(a_j)_{j=0}^{k-1}$ and $(q_i)_{i=0,-1}^{r-1}$ which remain $\mathbb{Q}$-linearly independent, no period vanishes along this path. }. While $A(t) > 0$, $F(X_t) = F(X_0) - \sqrt{-1}\alpha_0t$ must remain in $\bH$. But choosing $q_0$ arbitrarily small means that $t$ can be made large enough so that $\mathrm{Im}(F(t)) < 0$ and $A(t) >0$, a contradiction.
\end{proof}

\section{Proof of Theorems \ref{T:KraGenusZero} and \ref{T:KraGenusOne}: The converse in genus zero and one}\label{S:lowgenus}
\subsection{Theorem \ref{T:KraGenusZero} (Converse to Kra's theorem in genus zero)}

%This follows quickly from Proposition \ref{P:CylinderFreeGenus0}.

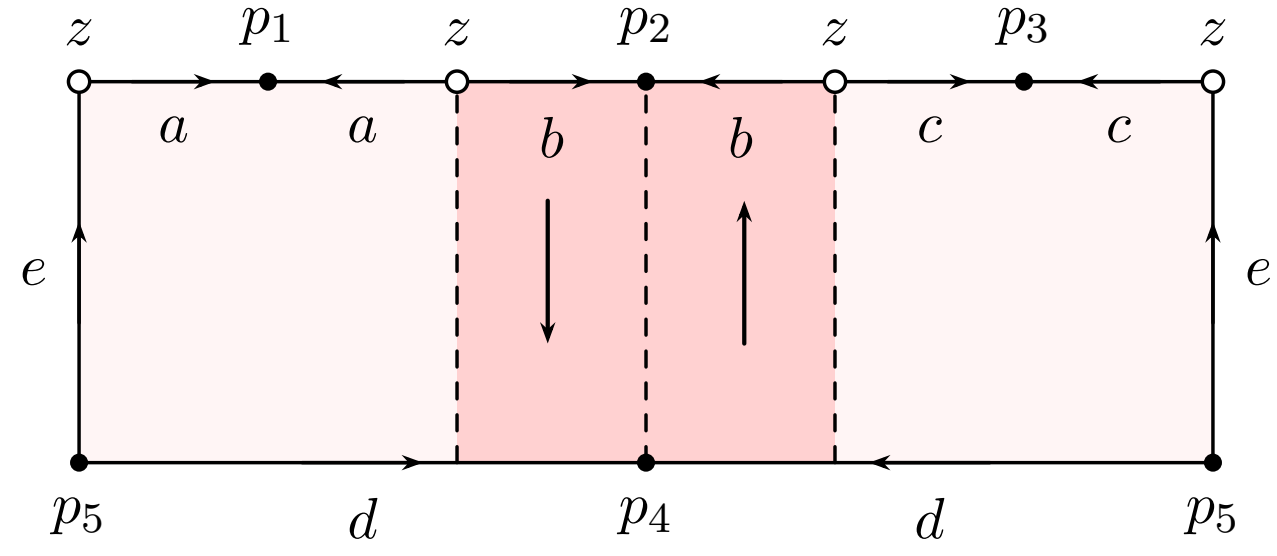
\begin{figure}[H]
\centering
\begin{tikzpicture}[
    x=1.25cm,y=1.05cm,
    draw=black,text=black,line width=.65pt,
    line cap=round,line join=round,
    >={Stealth[length=1.4mm,width=1mm]},
    every node/.style={font=\small}]
  % One horizontal cylinder: circumference 6, height 2.4.
  % Top: a a b b c c, paired by half-turns.
  % Bottom: d d, paired by a half-turn.
  % Vertical sides: e e, paired by translation.
  \fill[red!4] (0,0) rectangle (6,2.4);
  \fill[red!18] (2,0) rectangle (4,2.4);
  \draw (0,0) rectangle (6,2.4);
  \draw[dashed] (2,0)--(2,2.4);
  \draw[dashed] (3,0)--(3,2.4);
  \draw[dashed] (4,0)--(4,2.4);

  % Arrows and labels specify all edge identifications.
  \foreach \u/\lab in {0/a,2/b,4/c}{
    \draw[->] ({\u+.28},2.4)--({\u+.72},2.4);
    \draw[->] ({\u+1.72},2.4)--({\u+1.28},2.4);
    \node[below] at ({\u+.5},2.32) {$\lab$};
    \node[below] at ({\u+1.5},2.32) {$\lab$};
  }
  \draw[->] (1.18,0)--(1.82,0);
  \draw[->] (4.82,0)--(4.18,0);
  \node[below] at (1.5,-.08) {$d$};
  \node[below] at (4.5,-.08) {$d$};
  \draw[->] (0,.88)--(0,1.52);
  \draw[->] (6,.88)--(6,1.52);
  \node[left] at (-.05,1.2) {$e$};
  \node[right] at (6.05,1.2) {$e$};

  % The four open circles represent the same simple zero.
  \foreach \u in {0,2,4,6}{
    \draw[fill=white] (\u,2.4) circle (2pt);
    \node[above] at (\u,2.49) {$z$};
  }
  % Five simple poles; p_5 appears twice at the cut.
  \foreach \u/\lab in {1/1,3/2,5/3}{
    \fill (\u,2.4) circle (1.7pt);
    \node[above] at (\u,2.49) {$p_{\lab}$};
  }
  \foreach \u/\lab in {0/5,3/4,6/5}{
    \fill (\u,0) circle (1.7pt);
    \node[below] at (\u,-.08) {$p_{\lab}$};
  }
  % The two strips form one vertical cylinder.
  \draw[->,thick] (2.48,1.65)--(2.48,.75);
  \draw[->,thick] (3.52,.75)--(3.52,1.65);
\end{tikzpicture}
\caption{This is a surface in $\cQ(1, -1^5)$ with one horizontal cylinder. Letters indicate edge identification, the zero is marked with a $z$ and the poles by $p_i$. The darker shaded vertical cylinder corresponding to the saddle connection between $p_2$ and $p_4$ is the one created by the proof of Proposition \ref{P:1Cylinder}.}
\label{F:PTS}
\end{figure}

\begin{prop}\label{P:1Cylinder}
The only genus-zero cylinder-free Teichm\"uller curve is $\cQ(-1^4)$.

\end{prop}
\begin{proof}
A cylinder-free Veech surface has exactly one
cylinder in every cylinder direction. Suppose that such a genus-zero surface $X$ is
horizontally periodic. Then it necessarily has poles on opposite
boundaries. To see this, note that cutting along the horizontal cylinder core curve produces two spheres with boundary and hence the sum of orders of singularities on either cylinder boundary is $-2$. Shear the horizontal cylinder to line up two such poles
on opposite boundaries. This produces a vertical cylinder that will
not cover the surface unless it belongs to $\cQ(-1^4)$. See Figure \ref{F:PTS}.
\end{proof}

\begin{prop}\label{P:CylinderFreeGenus0}
The only cylinder-free invariant subvariety of genus-zero surfaces are strata.
\end{prop}
\begin{proof}
Let $\cM$ be such a cylinder-free invariant subvariety. By Apisa--Wright
\cite[Theorem~8.1(1)--(2)]{ApisaWrightGeminal}, $\cM$ is a full locus of degree $d_\cM$ covers of genus zero surfaces whose holonomy covers form an $h$-geminal subvariety. Moreover, the boundary components of $h$-geminal subvarieties are $h$-geminal and consist of connected surfaces. The degree\footnote{In principle, surfaces in $\cM$ may admit multiple half-translation covers of spheres. This degree refers to the degree of an ``optimal'' cover. For a definition see \cite{ApisaWrightGeminal}.} $d_\cM$ of the covers is also unchanged by generic cylinder collapses to the boundary. See \cite[Lemmas~3.11, 4.4, 4.6, 4.18 and Theorem~8.1(4)]{ApisaWrightGeminal}. Applying sufficiently many of these cylinder collapses causes $\cM$ to become a cylinder-free Teichm\"uller curve of genus-zero surfaces and hence $\cQ(-1^4)$ by Proposition \ref{P:1Cylinder}, so $d_\cM=1$ and so $\cM$ is a full stratum. 
\end{proof}

\begin{proof}[Proof of Theorem \ref{T:KraGenusZero}:]
This is immediate from Corollary \ref{C:CylinderFree},
Proposition \ref{P:CylinderFreeGenus0}, and Theorem \ref{T:Generic}.
\end{proof}

\subsection{Theorem \ref{T:KraGenusOne} (Converse to Kra's theorem in genus one)}
As in the previous subsection, it suffices to classify cylinder-free
subvarieties in genus one.
Recall that the \emph{ornithorynque} is the
pillowcase-tiled surface (PTS) that is a $\bZ/3$-cyclic cover of
$\bP^1$ branched over three of the four poles of the
pillowcase.

\begin{lem}\label{L:UniqueOneCylinderPTS}
The ornithorynque is the only cylinder-free genus one PTS
\end{lem}
\begin{proof}
Let $X$ be a genus one cylinder-free PTS and let $n$ be the degree of its cover to a pillowcase in $\cQ(-1^4)$. This surfaces is specified by four permutations $(a,b,c,d)$ in $\mathrm{Sym}(n)$ whose product is the identity and so that if $(a',b',c',d')$ is in its
braid group orbit then $a'b'$ is an $n$-cycle (see \cite[Section 8.4]{ApisaWrightGeminal}).
As in the proof of Proposition \ref{P:1Cylinder}, it cannot be the case that opposite sides of a cylinder in $X$ have poles. Because any two distinct poles of the pillowcase can be placed on opposite boundaries of a cylinder by choosing an appropriate periodic direction, it follows that exactly one element of $(a,b,c,d)$ has a fixed point (note that $X$ must also have poles, which correspond to such fixed points). Suppose that $d$ is this element. 

For any choice of two elements $(e,f)$ from $(a,b,c)$, there is a
cylinder direction on $X$ so that one boundary of the cylinder has
singularities that are exactly the preimages of
the branch values corresponding to $e$ and $f$.
Since $X$ has genus one, these singularities must be zeros whose orders sum to $2$ (since the sum of zeros is $-2\chi$ where $\chi$ is the Euler characteristic of the subsurface containing these preimages formed by cutting the cylinder core curve).
It follows that each element of $(a,b,c)$ has cycle type $(3,2^k)$
where $k=\frac{n-3}{2}$.
Letting $d_i$ be the lengths of the cycles of $d$
counted with multiplicity, we therefore have
$\sum_i(d_i-2)=-3$ and $\sum_i d_i=n$.
So $L:=\frac{n+3}{2}$ is the number of cycles of $d$.

Let $\deg(\lambda)$ be the degree of the Hodge bundle over a (suitable compactification of) the family determined by the PTS. We will see in a moment that we actually don't need to pass to any compactification. By Eskin--Kontsevich--Zorich
\cite[Corollary~8(b) and Section 3.3]{EKZbig},
\[
72\deg(\lambda)
=N\left(n-\sum_{i=1}^4\sum_j\frac{1}{c_{ij}}\right)+\frac{3N}{n}
\]
where $N$ is the index of the Veech group of $X$ in $\mathrm{PSL}(2, \bZ)$ and where $(c_{ij})$ is a list of cycle lengths for the $i$th permutation.
For a cylinder-free PTS, all cylinder core curves are
separating. In genus one this means that pinching such a core
curve pinches off a $\bP^1$, which, in particular, does not pass
to the boundary of $\cM_1$. So cylinder-free PTS in genus one are isotrivial families.
Since they project to a point in $\cM_1$, $\deg(\lambda)=0$ and so,
\[
\sum_{ij}\frac{1}{c_{ij}}=n+\frac{3}{n}.
\]
Combining this with the known cycle types of $(a,b,c)$, we have,
\[
\sum_{i=1}^L d_i=n
\quad\text{and}\quad
\sum_{i=1}^L\frac{1}{d_i}=\frac{n+5}{4}+\frac{3}{n}.
\]
Therefore, 
\[
\frac{n+5}{4}+\frac{3}{n}
=\sum_i\frac{1}{d_i}
\geq\sum_i\frac{3-d_i}{2}
=\frac{3L-n}{2}
=\frac{n+9}{4}.
\]
This holds if and only if $\frac{3}{n}\geq1$, so $n\leq3$ and
hence $n=3$ since $k\geq0$.
So $a,b,c$ are $3$-cycles in $\mathrm{Sym}(3)$
and $d$ is the identity, which is precisely the monodromy of the
ornithorynque.
\end{proof}

\begin{lem}\label{L:RankOne}
Aside from strata, the only cylinder-free invariant subvarieties in genus one strata that are loci of covers of $\cQ(0,-1^4)$ are ones where a marked point is added to the ornithorynque.
\end{lem}
\begin{proof}
Let $\cM$ be such a subvariety and let $d_\cM$ be the degree of the cover of surfaces in $\cQ(0, -1^4)$. (For technical reasons this is actually the degree of an ``optimal'' cover). Just like in the proof of Proposition \ref{P:CylinderFreeGenus0}, this degree is constant upon passing to the boundary of $\cM$. Passing to the boundary therefore does not cause the genus to drop by Proposition \ref{P:1Cylinder}. We can collide the zero with any of the four poles and still get a cylinder-free cover of the same degree. This shows,
by Lemma \ref{L:UniqueOneCylinderPTS}, that our monodromy is in $\mathrm{Sym}(3)$ and that at least two poles, say $p_1$ and $p_2$, have monodromy $\sigma:=(1,2,3)$. It follows, by considering collisions of the zero with $p_1$ or $p_2$, that the monodromy around the marked point is either $\sigma^{-1}$ or trivial. Thus, the monodromy vector is either
$(\sigma,\sigma,\sigma,\sigma,\sigma^{-1})$ or $(\sigma,\sigma,\sigma,1,1)$.
The first option is genus three, so we are left with the second option as desired.
\end{proof}

\begin{prop}\label{P:CylinderFreeGenus1}
There are no proper cylinder-free invariant subvarieties of rank at least two that consist of genus one surfaces. Therefore, up to marked points, the only such invariant subvariety is the ornithorynque.
\end{prop}
\begin{proof}
By inductive degeneration arguments as in Proposition \ref{P:CylinderFreeGenus0} and Lemma \ref{L:RankOne}, it suffices to show that there is no proper cylinder-free invariant subvarieties consisting of genus one surfaces that are covers of surfaces in $\cQ(1, -1^5)$. But for those covers, colliding the zero and a pole produces trivial local monodromy at the resulting marked point (by Lemma \ref{L:RankOne}) and hence shows that all poles have the same monodromy cycle type, which, after colliding a zero with a pole, contradicts Lemma \ref{L:RankOne}.
\end{proof}

\begin{proof}[Proof of Theorem \ref{T:KraGenusOne}:]
This is immediate from Corollary \ref{C:CylinderFree},
Proposition \ref{P:CylinderFreeGenus1}, and Theorem \ref{T:Generic}.
\end{proof}

\section{Proof of Theorem \ref{T:JSsurface}: Cylindrical stability and merging} \label{S:Construction}
%A \emph{Jenkins-Strebel differential} is a quadratic differential whose nonsingular horizontal trajectories are compact, i.e. one that is horizontally periodic. 

Fix a component $\cC$ of a stratum of quadratic differentials $\cQ(\kappa)$. A cylinder is \emph{simple} if each of its boundary components is a single saddle connection of multiplicity one. A cylinder is a \emph{simple envelope} if one of its boundary components is a single saddle connection of multiplicity one and the other is a single saddle connection of multiplicity two. See \cite[Section 4]{ApisaWrightDiamonds} for illustrations. 

A surface $(X,q) \in \cQ(\kappa)$ is \emph{cylindrically stable} if it realizes the maximum number possible of horizontal cylinders in the stratum, which is $g+ \frac{m_{odd}}{2} + m_{even}-1 $ where $m_{odd}$ and $m_{even}$ are the number of odd and even entries of $\kappa$ respectively. Cylindrically stable surfaces satisfies Properties (2)-(5) of Theorem \ref{T:JSsurface} (Property (4) is vacuous in this case). Property (2) follows from definition, Property (3) follows from Lemma \ref{L:PeriodCoord}, and Property (5) follows from the fact each $q_i$ is a horizontal core curve.

\begin{lem} \label{L:PeriodCoord}
For a cylindrically stable surface, the collection of all cross curves and rank many horizontal cylinder core curves forms a basis of period coordinates.
\end{lem}
\begin{proof}
It suffices to find $\dim(\cQ)$-many linearly independent periods. Consider first the $(\mathrm{rank}(\cQ) + \mathrm{rel}(\cQ))$ many cross curves of horizontal cylinders. They generate a linearly independent set in period coordinates since each one can be altered by a cylinder dilation while fixing all others. The periods of the horizontal core curves of a cylindrically stable surface span a $\mathrm{rank}(\cQ)$ dimensional subspace in period coordinates by \cite[Lemma 7.10(4)]{ApisaWrightHighRank}. They are linearly independent from the previously chosen periods again by considering cylinder dilations. We are done because $\dim(\cQ)= (2\mathrm{rank}(\cQ) + \mathrm{rel}(\cQ))$.
\end{proof}
Hence, it suffices to find a cylindrically stable surface which satisfies Property (1) of Theorem \ref{T:JSsurface} whenever possible. In the sequel, $z$ is the odd order zero that has been specified in the statement of Theorem \ref{T:JSsurface}.

 \begin{lem} \label{L:genericcyl}
 Theorem \ref{T:JSsurface} holds whenever $\cC$ contains a surface with a simple envelope or simple cylinder $C$ that has $z$ on one boundary component and a distinct odd order singularity on the other.
% Suppose that there is a surface in $\cC$ with a simple envelope or simple cylinder $C$ which has $z$ on one boundary component and a distinct odd order singularity on the other. Then Theorem \ref{T:JSsurface} holds for $\cC$.
\end{lem}
\begin{proof}
Rotate so that $C$ is horizontal and then perturb the surface in $\cC$ while keeping $C$ horizontal to pass to a horizontally periodic cylindrically stable surface by \cite[Lemma 7.10 (6)]{ApisaWrightHighRank}. Being a simple envelope or simple is a local condition, so $C$ remains a simple envelope or simple on this surface. Setting $C_0 := C$, we are done.
%
% Call that cylinder $C$. We can then rotate the surface so that $C$ is horizontal. In a small neighborhood of $\cC$, $C$ remains a simple cylinder or simple envelope with the same singularities on its boundary. We can find a surface in this neighborhood using \cite[Lemma 7.10 (6)]{ApisaWrightHighRank} that is cylindrically stable with $C$ a desired cylinder (up to shearing $C$ so that we have a vertical saddle connection between $z$ and an odd order singularity).
\end{proof}

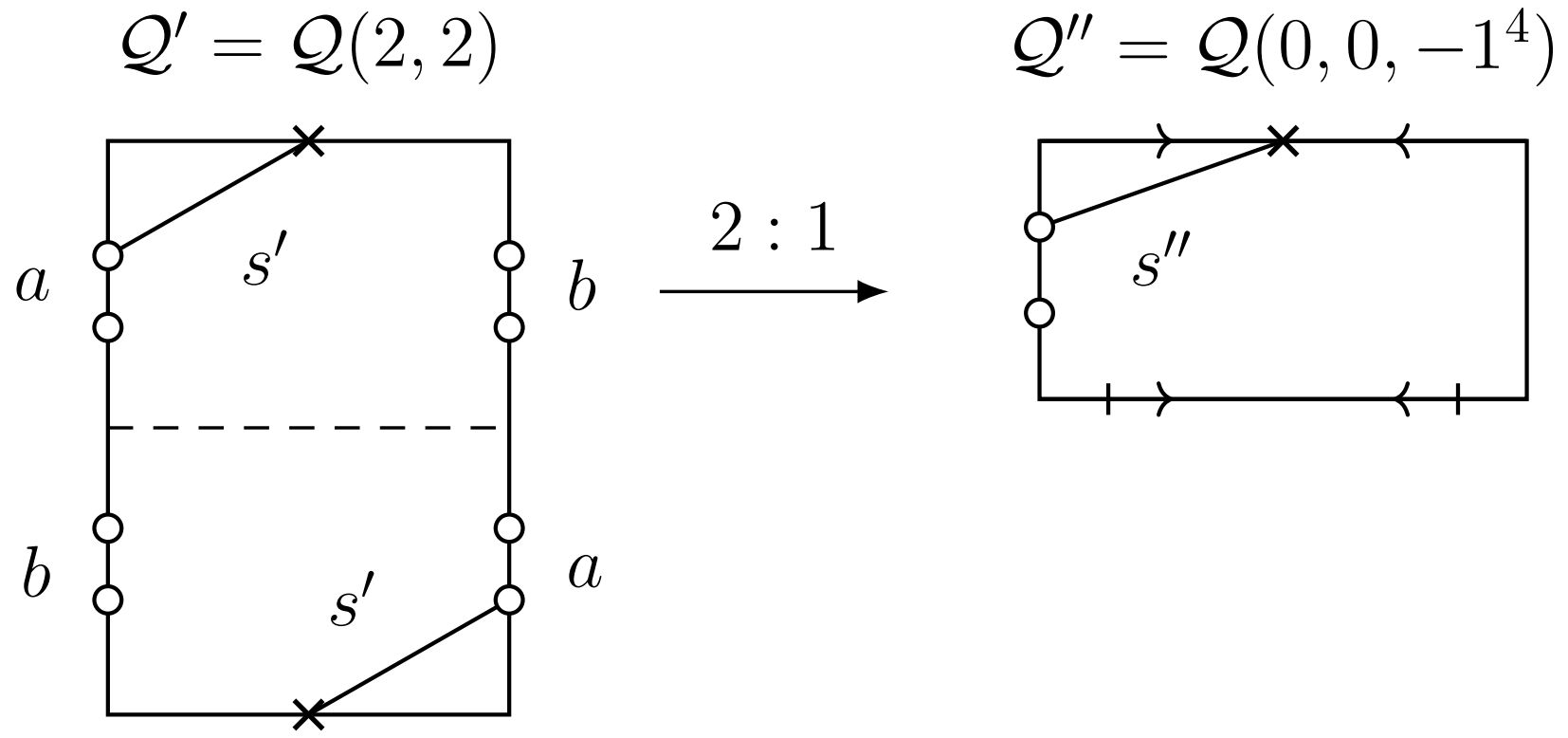
\begin{figure}
     \centering
\begin{tikzpicture}[x=1cm,y=1cm,line width=.8pt,
 every node/.style={font=\large},
 directed/.style={postaction={decorate},decoration={markings,
 mark=at position .55 with {\arrow{>}}}}]
% Hyperelliptic surface: paired vertical sides are labeled a and b.
\node at (1.4,4.65) {$\mathcal Q'=\mathcal Q(2,2)$};
\draw (0,0) rectangle (2.8,4);
\draw[dash pattern=on 5pt off 4pt,line width=.65pt] (0,2)--(2.8,2);
\node[left=7pt] at (0,3) {$a$};
\node[left=7pt] at (0,1) {$b$};
\node[right=7pt] at (2.8,3) {$b$};
\node[right=7pt] at (2.8,1) {$a$};
% The indicated saddle connections, both labeled sigma'.
\draw (0,3.2)--(1.4,4) node[midway,below right=2pt] {$s'$};
\draw (1.4,0)--(2.8,.8) node[midway,above left=2pt] {$s'$};
% Open-circle markings on the paired edges.
\foreach \x in {0,2.8}
 \foreach \y in {.8,1.3,2.7,3.2}
  {\draw[fill=white] (\x,\y) circle (2.7pt);}
% Cross markings for the two starred points.
\foreach \y in {0,4}{
 \draw[line width=1pt] (1.30,\y-.10)--(1.50,\y+.10);
 \draw[line width=1pt] (1.30,\y+.10)--(1.50,\y-.10);
}
% Degree-two quotient map.
\draw[-{Latex[length=2.2mm]},line width=.65pt] (3.85,2.95)--(5.45,2.95)
 node[midway,above=4pt] {$2:1$};
% Genus-zero quotient, including its two marked regular points.
\begin{scope}[shift={(6.5,2.2)}]
 \node at (1.7,2.45) {$\mathcal Q''=\mathcal Q(0,0,-1^{4})$};
 \draw (0,0) rectangle (3.4,1.8);
 % The upper and lower boundaries each fold onto themselves.
 \draw[directed] (0,1.8)--(1.7,1.8);
 \draw[directed] (3.4,1.8)--(1.7,1.8);
 \draw[directed] (0,0)--(1.7,0);
 \draw[directed] (3.4,0)--(1.7,0);
 \draw (.48,-.11)--(.48,.11);
 \draw (2.92,-.11)--(2.92,.11);
 \draw (0,1.2)--(1.7,1.8) node[midway,below=5pt] {$s''$};
 \foreach \y in {.6,1.2} \draw[fill=white] (0,\y) circle (2.7pt);
 \draw[line width=1pt] (1.60,1.70)--(1.80,1.90);
 \draw[line width=1pt] (1.60,1.90)--(1.80,1.70);
\end{scope}
\end{tikzpicture}
     \caption{The construction with surfaces in $\cQ(2,2)$ covering surfaces in $\cQ(0^2, -1^4)$.} 
     \label{F:Hyp}
 \end{figure}

\begin{lem} \label{L:hyperelliptic}
If $\cC$ is hyperelliptic, then Theorem \ref{T:JSsurface} holds. In particular, we can find a simple cylinder with $z$ and its (distinct) image under the hyperelliptic involution on its boundary.
\end{lem}
\begin{proof}
Since $z$ is of odd order, it is distinct from its image under the hyperelliptic involution and $\cQ = \cQ(z, z, \kappa')$ where $\kappa' = (2a+2), (a,a)$ for some integer $a \geq -1$. Consider the hyperelliptic component of $\cQ' := \cQ(2z, \kappa')$, which is a full locus of covers of $\cQ'' := \cQ(z-1, a, -1^{3+z+a})$. Find a surface $X$ in $\cQ''$ with a saddle connection $s''$ joining $(z-1)$ to one of the $3+z+a$ poles. (See Figure \ref{F:Hyp}). On the corresponding double cover in $\cQ'$, this becomes a saddle connection $s'$ that joins the zero of order $2z$ to itself and that is fixed by the hyperelliptic involution. (Such a surface with this saddle connection exists either by merging arguments relying on \cite{BCGGM} as in \cite{ApisaAygun} or by the solution of the finite blocking problem for strata in \cite{ApisaWright}). Gluing in a cylinder to the saddle connection $s'$ produces a surface in a hyperelliptic component of $\cQ(z,z, \kappa')$ (see \cite[Lemma 6.3]{ApisaWright}) with a simple cylinder with the two desired zeros on its boundary.
%
% Let $J$ be the hyperelliptic involution. Since $z$ is of odd order, $z$ and $J(z)$ are distinct. Merge $z$ with $J(z)$ to form a singularity $x$ and glue a cylinder into a closed saddle connection based at $x$ and fixed by the hyperelliptic involution. We can find this fixed saddle connection by lifting a saddle connection between a simple pole and the projection of $x$ in the quotient by $J$. The cylinder we glue in is our desired cylinder. The claim then follows from Lemma \ref{L:genericcyl}.
\end{proof}

In the sequel we will use Lemma \ref{L:hyperelliptic} to assume that $\cC$ is not hyperelliptic. Since marked regular points in $\cC$ can be freely moved, we may safely forget them. 

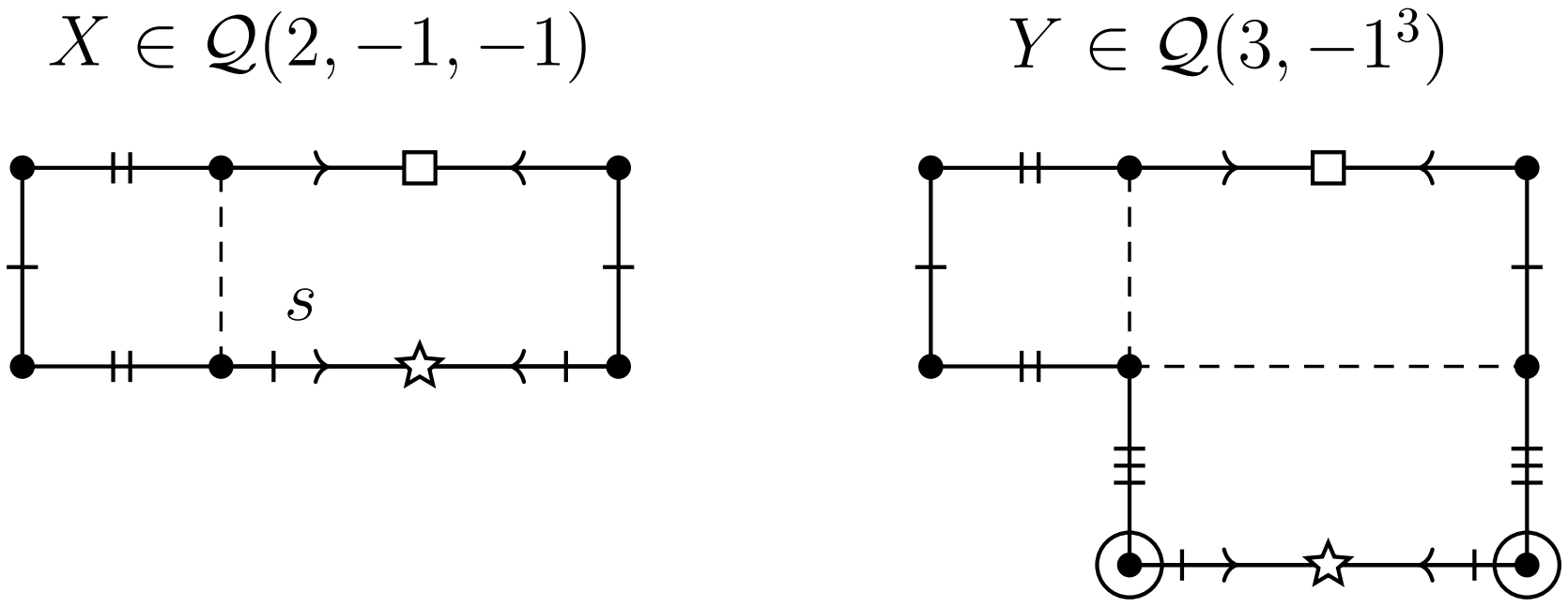
\begin{figure}
    \centering
\begin{tikzpicture}[
  x=1cm,y=1cm,line width=.8pt,
  every node/.style={font=\large},
  directed/.style={
    postaction={decorate},
    decoration={
      markings,
      mark=at position .55 with {\arrow{>}}
    }
  },
  starpoint/.style={
    star,star points=5,star point ratio=2.2,
    draw,fill=white,minimum size=9pt,inner sep=0pt
  }
]

% Left surface.
\node at (2.1,2.25) {$X\in\mathcal Q(2,-1,-1)$};

\draw (0,0) rectangle (4.2,1.4);
\draw[dash pattern=on 4pt off 3pt,line width=.65pt]
  (1.4,0)--(1.4,1.4);

% Folded horizontal edges.
\draw[directed] (1.4,1.4)--(2.8,1.4);
\draw[directed] (4.2,1.4)--(2.8,1.4);
\draw[directed] (1.4,0)--(2.8,0);
\draw[directed] (4.2,0)--(2.8,0);

% Single ticks on the paired vertical sides.
\draw (-.11,.7)--(.11,.7);
\draw (4.09,.7)--(4.31,.7);

% Double ticks on the left horizontal pair.
\foreach \x in {.64,.76}{
  \draw (\x,-.11)--(\x,.11);
  \draw (\x,1.29)--(\x,1.51);
}

% Ticks on the folded lower boundary.
\draw (1.77,-.11)--(1.77,.11);
\draw (3.83,-.11)--(3.83,.11);
\node[above=5pt] at (1.96,0) {$s$};

% All ordinary vertices are solid dots.
\foreach \p in {
  (0,0),(0,1.4),
  (1.4,0),(1.4,1.4),
  (4.2,0),(4.2,1.4)
}{
  \fill \p circle (2.5pt);
}

% Square and star.
\draw[fill=white] (2.69,1.29) rectangle (2.91,1.51);
\node[starpoint] at (2.8,0) {};

% Right surface.
\begin{scope}[shift={(6.4,0)}]
  \node at (2.1,2.25) {$Y\in\mathcal Q(3,-1^{3})$};

  \draw
    (0,0)--(0,1.4)--(4.2,1.4)--
    (4.2,-1.4)--(1.4,-1.4)--(1.4,0)--cycle;

  % Dashed boundary of the simple envelope.
  \draw[dash pattern=on 4pt off 3pt,line width=.65pt]
    (1.4,0)--(4.2,0);

  % Dashed vertical segment.
  \draw[dash pattern=on 4pt off 3pt,line width=.65pt]
    (1.4,0)--(1.4,1.4);

  % Folded horizontal edges.
  \draw[directed] (1.4,1.4)--(2.8,1.4);
  \draw[directed] (4.2,1.4)--(2.8,1.4);
  \draw[directed] (1.4,-1.4)--(2.8,-1.4);
  \draw[directed] (4.2,-1.4)--(2.8,-1.4);

  % Single ticks on upper vertical sides.
  \draw (-.11,.7)--(.11,.7);
  \draw (4.09,.7)--(4.31,.7);

  % Double ticks on the left horizontal pair.
  \foreach \x in {.64,.76}{
    \draw (\x,-.11)--(\x,.11);
    \draw (\x,1.29)--(\x,1.51);
  }

  % Triple ticks on the envelope's vertical sides.
  \foreach \y in {-.58,-.7,-.82}{
    \draw (1.29,\y)--(1.51,\y);
    \draw (4.09,\y)--(4.31,\y);
  }

  % Ticks on the folded lower boundary.
  \draw (1.77,-1.51)--(1.77,-1.29);
  \draw (3.83,-1.51)--(3.83,-1.29);

  % Solid dots at all ordinary vertices.
  \foreach \p in {
    (0,0),(0,1.4),
    (1.4,0),(1.4,1.4),
    (4.2,0),(4.2,1.4),
    (1.4,-1.4),(4.2,-1.4)
  }{
    \fill \p circle (2.5pt);
  }

  % Circles around the two lower dots.
  \foreach \x in {1.4,4.2}{
    \draw (\x,-1.4) circle (6.5pt);
  }

  % Square and star.
  \draw[fill=white] (2.69,1.29) rectangle (2.91,1.51);
  \node[starpoint] at (2.8,-1.4) {};
\end{scope}

\end{tikzpicture}
     \caption{The envelope construction with surfaces in $\cQ(2,-1^2)$ being used to construct a surface with a simple envelope in $\cQ(3, -1^3)$.} 
     \label{F:SimpleEnvelope}
 \end{figure}

\begin{lem}\label{L:TwoPoles}
If $\kappa$ contains at least two poles, then there is a surface in $\cQ$ where $z$ and a pole lie on opposite boundaries of a simple envelope; so Theorem \ref{T:JSsurface} holds.
\end{lem}
\begin{proof}
Let $\cQ := \cQ(z, -1, -1, \kappa')$ be the original stratum. If $\kappa' = (z)$, then $\cQ$ has a hyperelliptic component and, in all cases (except for $\cQ = \cQ(1^2, -1^2)$, which is excluded by being hyperelliptic), it has exactly one nonhyperelliptic component by Lanneau \cite{Lconn}. The stratum $\cQ(z-1, -1, \kappa')$ is nonempty by Lanneau \cite{Lconn}. Start with a surface $X$ in $\cQ(z-1, -1, \kappa')$ that has a saddle connection $s$ joining $z-1$ to $-1$. (Such a surface with this saddle connection exists either by merging arguments relying on \cite{BCGGM} as in \cite{ApisaAygun} or by the solution of the finite blocking problem for strata in \cite{ApisaWright}; see Figure \ref{F:SimpleEnvelope}). Glue a simple envelope into $s$ to produce a surface $Y \in \cQ(z, -1, -1, \kappa')$. If $X$ belongs to a hyperelliptic component, then, since we are gluing in a simple envelope to a saddle connection that is not fixed by the involution, $Y$ does not belong to a hyperelliptic component. We are done by Lemma \ref{L:genericcyl}. 
\end{proof}

% \begin{figure}
%     \centering

%      \caption{The construction with surfaces in $\cQ(2,2)$ covering surfaces in $\cQ(0^2, -1^4)$.} 
%      \label{F:Hyp}
%  \end{figure} 

\begin{figure}
    \centering
\begin{tikzpicture}[line width=.8pt,
 every node/.style={font=\large},
 directed/.style={postaction={decorate},decoration={markings,
 mark=at position .55 with {\arrow{>}}}},
 starpoint/.style={star,star points=5,star point ratio=2.2,
 draw,fill=white,minimum size=9pt,inner sep=0pt}]
% Common outlines and edge identifications.
\foreach \offset in {0,6}{
 \begin{scope}[shift={(\offset,0)}]
  \draw (0,0)--(0,1.1)--(4.2,1.1)--(4.2,-2)
    --(1.4,-2)--(1.4,0)--cycle;
  \draw[dash pattern=on 4pt off 3pt,line width=.65pt]
    (1.4,1.1)--(1.4,0)--(4.2,0);
  \draw[directed] (1.4,1.1)--(2.8,1.1);
  \draw[directed] (4.2,1.1)--(2.8,1.1);
  \draw[directed] (1.4,-2)--(2.8,-2);
  \draw[directed] (4.2,-2)--(2.8,-2);
  \draw (-.11,.55)--(.11,.55);
  \foreach \x in {.64,.76}{
   \draw (\x,-.11)--(\x,.11);
   \draw (\x,.99)--(\x,1.21);
  }
  \draw (1.77,-2.11)--(1.77,-1.89);
  \draw (3.83,-2.11)--(3.83,-1.89);
  \foreach \p in {(0,0),(0,1.1),(1.4,0),(1.4,1.1),
                   (4.2,0),(4.2,1.1)}
   \fill \p circle (2.3pt);
  \draw[fill=white] (2.69,.99) rectangle (2.91,1.21);
  \node[starpoint] at (2.8,-2) {};
 \end{scope}
}
% Left surface: the indicated saddle connection s.
\draw (4.09,.55)--(4.31,.55);
\draw[line width=1.3pt] (1.4,-1)--(1.4,-2);
\node[left=6pt] at (1.4,-1.5) {$s$};
\foreach \x in {1.4,4.2}{
 \draw[fill=white] (\x,-1) circle (3pt);
 \fill (\x,-2) circle (2.1pt);
 \draw (\x,-2) circle (4pt);
}
% Right surface: the cut edges b and a.
\begin{scope}[shift={(6,0)}]
 \node[right=6pt] at (1.4,-1.5) {$b$};
 \node[left=6pt] at (4.2,-1.5) {$a$};
 \foreach \x in {1.4,4.2}{
  \foreach \y in {-1,-2}{
   \fill (\x,\y) circle (2.3pt);
   \draw (\x,\y) circle (6pt);
  }
 }
\end{scope}
% Separate rectangle with the matching edges.
\begin{scope}[shift={(11.5,-2)}]
 \draw (0,0) rectangle (1.1,1);
 \node[left=5pt] at (0,.5) {$a$};
 \node[right=5pt] at (1.1,.5) {$b$};
 \foreach \x in {.43,.55,.67}{
  \draw (\x,-.10)--(\x,.10);
  \draw (\x,.90)--(\x,1.10);
 }
\end{scope}
\end{tikzpicture}

     \caption{The slit construction using the slit $s$ based at the pole represented by the circled solid dot.} 
     \label{F:add in them slits!}
 \end{figure}
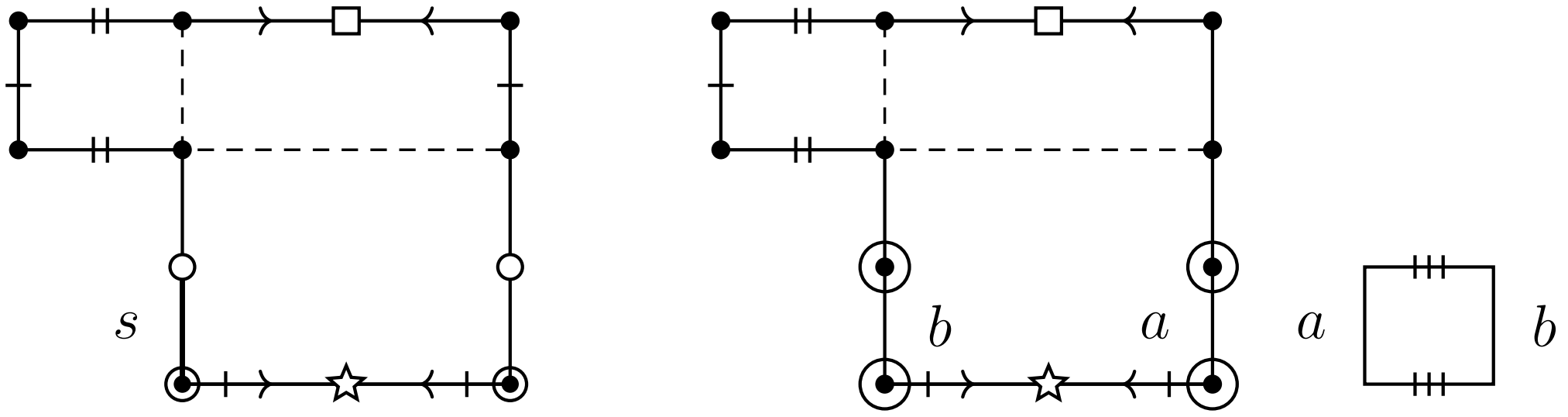

By Lemma \ref{L:TwoPoles}, we can from now on assume that there is at most one pole in $\kappa$. The \emph{(first) slit construction} is a way of raising the order of a singularity $p$ of order $k$ in a quadratic differential to order $k+4$ by slitting a line segment joining $p$ to a marked point and gluing in a flat cylinder. See Figure \ref{F:add in them slits!}.

\begin{lem}\label{L:SlitConstructionIsNice}
If $X$ is a cylindrically stable surface and $p$ is a singularity, then applying the slit construction to a vertical segment based at $p$ whose other endpoint does not lie on a horizontal saddle connection is cylindrically stable and belongs to a nonhyperelliptic component. 

Moreover, if $X$ contains a horizontal cylinder $C$ with distinct singularities $q$ and $r$ on opposite boundaries that may be merged by collapsing $C$, then the slit based at $p$ may be chosen so that the same holds after applying the slit construction.
\end{lem}
\begin{proof}
Recall that the maximum number of cylinder in a stratum of genus $g$ surfaces with $m_{odd}$ and $m_{even}$ singularities of odd and even order is $g + \frac{m_{odd}}{2} + m_{even}-1$. The slit construction doesn't change the number of even and odd order singularities, but it does increase genus by one. The slit construction applied to a vertical segment whose endpoint does not lie on a horizontal saddle connection can be seen to increase the number of horizontal cylinders by exactly one and hence preserves cylindrical stability (see Figure \ref{F:add in them slits!}). Nonhyperellipticity is preserved since the boundary of hyperelliptic components must be hyperelliptic and reversing the slit construction passes to a stratum with a marked point, which cannot be.

The second claim is clear when $p$ does not lie on the boundary of $C$ by choosing a sufficiently small slit. If $p$ does lie on the boundary of $C$ then choosing the slit to be contained in $C$ produces the desired result.
\end{proof}

% Moreover, using a sufficiently small vertical slit, it is always possible to preserve cylindrical stability. 

% By Lemma \ref{L:TwoPoles} there are either two odd order zeros $z,w \geq 0$ or the unique subset of odd order singularities of $\kappa$ is $(z, -1)$, in which case set $w = -1$. By Lemma \ref{L:1Mod4Case}, the remaining cases to consider are where (1) $z,w \geq 0$ and $z = w = 3$ mod $4$ or (2) $(z, -1)$ are the only odd order singularities in $\kappa$.
%
%Handle first: w = 1 mod 4 or w not a pole and z = 1 mod 4
%
%Remaining cases: w = 3 mod 4 and (w a pole or z = 3 mod 4), i.e. 
%Case 1: (z, -1)
%case 2: w = z = 3 mod 4. 

\begin{lem}\label{L:Two3Mod4Zeros}
If $\kappa$ contains two zeros, including $z$, of order congruent to $3$ mod $4$ and at most one pole, then there is a surface in a nonhyperelliptic component of $\cQ$ where Theorem \ref{T:JSsurface} holds.
\end{lem}
\begin{proof}
Set $z = 4k_1-1$ and $w = 4k_2 -1$ and let $\kappa'$ be the remaining elements of $\kappa$. Note that $\cQ' := \cQ(-1, -1, \kappa')$ is a nonempty stratum by Lanneau \cite{Lconn}. It has a nonhyperelliptic component precisely when $\kappa' \ne (2), (1,1)$. Take the simple envelope in $\cQ'$ produced by Lemma \ref{L:TwoPoles} when $\cQ'$ is nonhyperelliptic and take a complex envelope in $\cQ(2, -1, -1)$ and $\cQ(1,1,-1,-1)$ on a cylindrically stable surface otherwise. Apply the slit construction $k_1$ times to one pole and $k_2$ to the other to conclude by Lemma \ref{L:SlitConstructionIsNice}. 
\end{proof}

\begin{lem}\label{L:Two3Mod4Zeros2}
The previous lemma holds in all components of $\cQ$.
\end{lem}
\begin{proof}
By Lanneau \cite{Lconn} we need only consider $\cQ = \cQ(6,3,3), \cQ(3^4),\cQ(3^3, -1)$ which all have two nonhyperelliptic components. Let $z$ and $w$ be the two zeros congruent to $3$ mod $4$. Let $u$ be a third zero, which we can split to form $(u-3,-1)$, thereby passing to a boundary stratum $\cQ' = \cQ(3,-1,3,3), \cQ(0,-1,3^3), \cQ(0, -1, 3,3, -1)$ by \cite[Proposition 10.4]{ApisaAygun}. By Lemmas \ref{L:Two3Mod4Zeros} and \ref{L:hyperelliptic}, there is a simple cylinder between $z$ and $w$ when $\cQ' = \cQ(0, -1, 3, 3, -1)$. This remains a simple cylinder upon smoothing back into the original stratum by \cite[Lemma 2.15(1)]{MirWri}. Now that we have the claim for $\cQ' = \cQ(3^3, -1)$, the remaining two strata follow analogously. 
\end{proof}

% \begin{lem}\label{L:Certificate}
% The conclusion of Theorem \ref{T:JSsurface} holds if we can find a surface in $\cQ$ with a cylinder $C$ with a pole $p$ and one boundary that has a $z-z$ saddle connection $s_1$ and a $z-p$ saddle connection $s_2$ and where the opposite boundary contains $s_1$ and another $z-z$ saddle connection $s_3$ that is homotopic to the core curve of a cylinder $D$
% \end{lem}
% \begin{proof}
% By \cite[Lemma 7.10]{ApisaWrightHighRank}, perturb the surface so that $D$ and $C$ are kept horizontal and so that we pass to a cylindrically stable surface. A priori, this might involve altering the saddle connections on the boundary of $C$. However, $s_3$ is unaltered as it is homotopic to the core cure of $D$. Therefore, $s_2$ is unaltered since $s_2 + s_3$ is homologous to the core curve of $C$. The same reasoning shows that $s_1$ is unaltered.
% \end{proof}

\begin{figure}
     \centering
     \includegraphics[width =.65\textwidth]{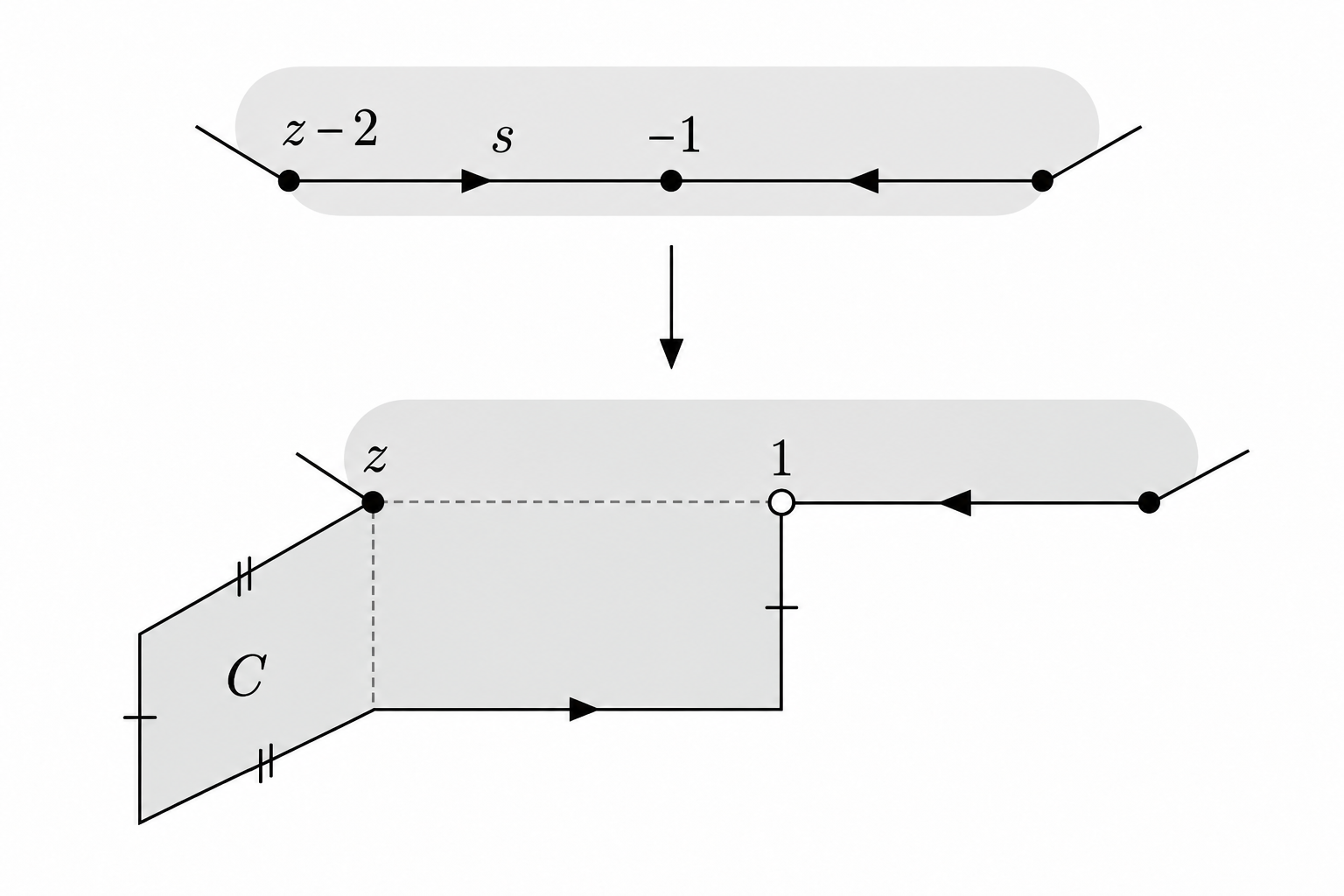}
     \caption{Creating a zero of order $1$.} 
     \label{F:Glue in subsurface}
 \end{figure}

\begin{lem}\label{L:1Mod4Case}
Suppose that $\cQ$ contains two zeros $z,w$, one of which is congruent to $1 \mod 4$ and where $\kappa \ne (a, 3)$ for some $a \in \bZ$. Then there is a surface in a nonhyperelliptic component of $\cQ$ with a simple cylinder $C$ with $z$ and $w$ on its boundary. 
\end{lem}
\begin{proof}
Without loss of generality, suppose that $ w \equiv 1 \mod 4$, $z \neq 3$. Write $\cQ = \cQ(z,w, \kappa')$. Let $X$ be a surface in $\cQ' := \cQ(z-2, -1, \kappa')$ with a saddle connection $s$ between $z-2$ and $-1$. Note that $\cQ$ is nonempty unless $\kappa = (a,3)$ for some $a \in \bZ$. The surface with the desired saddle connection can be found by the solution to the finite blocking problem for strata in \cite{ApisaWright}. Glue in a slit torus to $s$ as in Figure \ref{F:Glue in subsurface} to form a surface $Y \in \cQ(z, 1, \kappa')$ that contains a simple cylinder $C$ with $z$ and $1$ on its boundary. If $w = 4n+1$, then apply the slit construction to $1$, $n$ times with slits that do not enter $C$. We are done by Lemma \ref{L:genericcyl}.
 % We obtain a surface in the stratum $\mathcal{Q}(z,1,\underline{\kappa})$ with a cylinder $C$ that is simple and between $z$ and $1$. Then take $n$ marked points near $1$ away from $C$ where $w = 4n+1$, and glue in $n$ slit tori into the short saddle connections between the marked points and $1$. We constructed a surface in $\mathcal{Q}(\kappa)$ with our $C$ becoming a simple cylinder between $z$ and $w$. %See Figure \ref{F:add in them slits!}. %\begin{figure}
  %   \centering
   %  \includegraphics[width =.65\textwidth]{preArrival notes/Topology Surface Diagram with Insets.png}
    % \caption{The slit construction.} 
     %\label{F:add in them slits!}
 %\end{figure} Moreover, unless $\kappa= (9,3)$, $\mathcal{Q}(\kappa)$ is connected under our assumptions so this surface lies in $\cC$.
\end{proof}

%Special cases: 
%
% g=3: Q(9, -1); Q(6,3,-1); Q(3,3,3,-1). Only first is a problem. 
%
% g=4: Q(12); Q(9,3); Q(6,6); Q(6,3,3); Q(3,3,3,3). Only second is a problem.
%

\begin{lem} \label{L:simplezeroandpole}
    Theorem \ref{T:JSsurface} holds for the strata $\cQ = \cQ(1, -1, \kappa')$ where $\kappa'\neq (4),(2,2)$ consists of a set of even order zeros.
\end{lem}
\begin{proof} 
%
% Juliet's proof 1
%
%Merge the pole with an element of $\kappa'$ to create an odd order zero $w$. We now are in the case of Lemma \ref{L:1Mod4Case}. Consider a cylindrically stable surface $Y$ on this boundary stratum with a horizontal simple cylinder $C$ between $z$ and $w$. 
Since $\cQ(1,-1)$ is empty, assume that $\kappa'$ is nonempty. Write $\kappa' = (z') \cup \kappa''$. Set $w = z'-1$. By Lemma \ref{L:1Mod4Case}, let $Y$ be a cylindrically stable surface in $\cQ(w,1,\kappa'')$ with a horizontal simple cylinder $C$ with $z$ on one boundary and $w$ on the other. 

There is a local operation that replaces odd order zeros $w$ with a simple pole and a zero of order $w+1$. Applied to periodic surfaces, a cylinder may be specified on whose boundary the simple pole appears. All other cylinders persist. The construction is recalled in Figure \ref{F:OddToEvenAndPole}.

Let $X$ be the result of applying the operation to $Y$ so that the pole appears on the boundary of $C$. The new cylinder $C_X$ has one boundary, containing the simple zero, consisting of a single saddle connection, and the other, consisting of two. Call these two saddle connections $s$, which joins $z'$ to itself and $s'$, which joins $z'$ to the pole. Collapsing $C_X$ in the direction of a cross curve connecting the simple zero and pole merges them. 

Note that the cylinder cross curves and core curves form a basis of period coordinates on $Y$ by Lemma \ref{L:PeriodCoord}. So this basis extends to a basis for $X$ by adding in the vanishing cycle, $s'$, of the reverse operation that takes $X$ to $Y$. Since $s+s'$ is the core curve of $C_X$, the basis on $Y$ can also be extended by adding in $s$, which does not intersect the vanishing cycle $t$ of the cylinder collapse of $C_X$. It is straightforward to check that the polyplane around $X$ is determined by fixing the periods of $s$ and the cylinder core curves. This verifies all conditions in Theorem \ref{T:JSsurface}.
%
% Juliet's proof 2
%
% We can now open up the node at $w$ so that $C$ has a single saddle connection based at $z$ on one boundary component and $w+1$ and $-1$ on the other. In particular, the latter boundary component has a short saddle connection $s$ between $w+1$ and $-1$ and a closed saddle connection, we call $s'$, based at $w+1$. Collapsing $C$ in the direction of a cross curve between $z$ and the pole merges $a$ and $-1$. The surface remains horizontally periodic because smoothing consisted of introducing a short horizontal saddle connection and perturbing the periods of other horizontal saddle connections along $\mathbb{R}$. Hence, we create a horizontally periodic surface with $C$ as our desired cylinder. 
%
%%%%%%%%%%%%%%%
%
\begin{figure}
    \centering
\begin{tikzpicture}[
  scale=0.85,
  line cap=round,
  line join=round,
  every node/.style={font=\scriptsize},
  boundary/.style={draw=black!55,thick},
  seam/.style={draw=black,thick},
  zero/.style={circle,fill=blue!70!black,inner sep=1.5pt},
  pole/.style={circle,fill=orange!90!black,inner sep=1.5pt}
]

\def\R{1}
\def\e{0.4}

% Horizontal edge with tick marks and a length label.
\def\edge#1#2#3#4{%
  \draw[seam] (#1,0)--(#2,0)
    node[midway,below=5pt] {$#4$};
  \pgfmathsetmacro{\midpoint}{(#1+#2)/2}
  \foreach \j in {1,...,#3}{%
    \pgfmathsetmacro{\tickx}{\midpoint+0.055*(\j-(#3+1)/2)}
    \draw[seam]
      (\tickx-0.02,-0.055)--(\tickx+0.02,0.055);
  }
}

\def\capboundary{%
  \draw[boundary] (\R,0)
    arc[start angle=0,end angle=180,radius=\R];
}

% Before: three half-circles.
\begin{scope}[shift={(0,0)}]
  \capboundary
  \edge{-\R}{0}{1}{1}
  \edge{0}{\R}{2}{1}
  \node[zero] at (0,0) {};
\end{scope}

\begin{scope}[shift={(2.3,0)}]
  \capboundary
  \edge{-\R}{0}{2}{1}
  \edge{0}{\R}{3}{1}
  \node[zero] at (0,0) {};
\end{scope}

\begin{scope}[shift={(4.6,0)}]
  \capboundary
  \edge{-\R}{0}{3}{1}
  \edge{0}{\R}{1}{1}
  \node[zero] at (0,0) {};
\end{scope}

% Construction arrow.
\draw[->,thick] (5.9,0.45)--(6.9,0.45);

% After: three modified half-circles.
\begin{scope}[shift={(8.2,0)}]
  \capboundary
  \edge{-\R}{-\e}{1}{1-\epsilon}
  \edge{-\e}{\R}{2}{1+\epsilon}
  \node[zero] at (-\e,0) {};
\end{scope}

\begin{scope}[shift={(10.5,0)}]
  \capboundary
  \edge{-\R}{\e}{2}{1+\epsilon}
  \edge{\e}{\R}{3}{1-\epsilon}
  \node[zero] at (\e,0) {};
\end{scope}

\begin{scope}[shift={(12.8,0)}]
  \capboundary
  \edge{-\R}{-\e}{3}{1-\epsilon}
  \edge{\e}{\R}{1}{1-\epsilon}

  % Two epsilon edges.
  \draw[seam] (-\e,0)--(0,0)
    node[midway,below=5pt] {$\epsilon$};
  \draw[seam] (\e,0)--(0,0)
    node[midway,below=5pt] {$\epsilon$};

  % Small interior arrowheads pointing towards the pole.
  \draw[seam,-{Stealth[length=2.2pt,width=2pt]}]
    (-0.32,0)--(-0.16,0);
  \draw[seam,-{Stealth[length=2.2pt,width=2pt]}]
    (0.32,0)--(0.16,0);

  \node[zero] at (-\e,0) {};
  \node[zero] at (\e,0) {};
  \node[pole] at (0,0) {};
\end{scope}
\end{tikzpicture}
    \caption{Converting a simple odd-order (blue) zero of order $k$ to a (blue) zero of order $k+1$ and an (orange) simple pole. Edge labels indicate lengths.}
    \label{F:OddToEvenAndPole}
\end{figure}
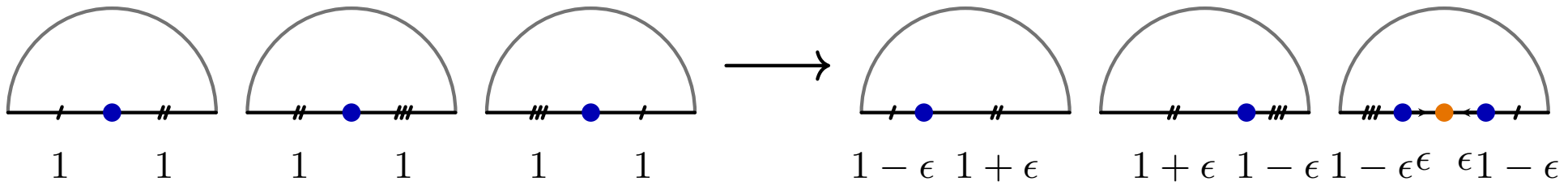
%%%%%%%%%%%%%%%%
%
% Juliet's Proof Part 3
%
%
% The last step is to verify Conditions (4) and (5) of Theorem \ref{T:JSsurface}. For Condition (4), we can take $q_{-1}$ to be $s'$ which is disjoint from $a_0$ and does not share endpoints. In particular, because the $\cup_i a_i$ and $\cup_{i>-1} q_i$ form period coordinates on the boundary by collapsing $s$, period coordinates for $\cQ$ can be given by $\cup_i a_i$, $\cup_{i>-1} q_i$, and $s$. Since $s,s',q_0$ bound an annulus and $q_0$ is a basis element, we can replace $s$ with $s'$ in our coordinate system. For Condition (5), note that locally varying any $a_i$ corresponds to a cylinder deformations about the horizontal cylinders, so fixing only each $q_i$ preserves horizontal periodicity.
\end{proof}

\begin{rem}
Lemma \ref{L:simplezeroandpole} is necessary since in the strata $\cQ(1,-1,\kappa)$ where $\kappa$ consists of even order zeros, there is no cylindrically stable surface that contains the pole on one boundary and the simple zero on the other. 
\end{rem}

\begin{proof}[Proof of Theorem \ref{T:JSsurface}:]  By Lemmas \ref{L:TwoPoles}, \ref{L:Two3Mod4Zeros2},  \ref{L:1Mod4Case}, and \ref{L:simplezeroandpole}, it remains to consider strata of the form $\cQ(w, 3)$ and of the form $\cQ = \cQ(z, -1, \kappa')$ where $\kappa'$ consists of (a possibly empty set of) even order zeros and $z >1$.

%nd the strata of the form $\cQ(w,3)$. We will handle the former case first. Write $\kappa' = (4k_i)_{i=1}^m \cup (4\ell_j+2)_{j=1}^n$. 
We will first show that it suffices to prove the claim for $\cQ(5,4,-1), \cQ(5,2,2,-1)$, $\cQ(9,-1)$, $\cQ(3,2,-1)$, $\cQ(5,3)$, and both components of $\cQ(9,3)$. Suppose the claim held in these strata. By Lemma \ref{L:SlitConstructionIsNice}, the slit construction give us the claim in strata where we have added multiples of $4$ to the order of any singularity. By applying the slit construction to a marked regular point, this includes creating zeros whose orders are multiples of $4$. Likewise, pairs of double zeros can be added by slitting a surface along a line segment and gluing this slit into a slit in a flat torus. Their orders can be raised by multiples of $4$ by the first slit construction. The second slit construction satisfies Lemma \ref{L:SlitConstructionIsNice} using a nearly identical proof. This reduces our claim to the strata listed above when $\cQ$ has a unique nonhyperelliptic component. However, the only strata with non-unique nonhyperelliptic components that we have left to analyze are $\cQ(6,3,-1)$ and $\cQ(9,-1)$. These have a unique component where the zero and pole can be merged. 

% Once we handle Q(3,2,-1) we may suppose that z = 1 mod 4
% Once we handle Q(9,-1) we may suppose that z = 1, 5. 
% We may suppose that kappa' is either (empty), (4), (2,2)
% Q(5,-1) is not possible
% Q(5,4,-1) is possible
% Q(5,2,2,-1) is possible and should be drawn in such a way that a rel deformation gives the pic of the previous. We could also do Q(5, -1,-1,-1)
%
% Now we are left with Q(1,-1,kappa)
% But we know that we are not genus 2, so the base cases are 
% kappa = (4,2,2), (4,4), (2,2,2,2)
%
% 
%
%
%To see this we may add multiples of four to the order of any zero using the first slit construction. We may add a zero of order $4$ by applying the slit construction to a marked point. 

To handle the case of $\cQ(5,3)$, apply the slit construction to a surface in $\cQ(5,-1)$ as in Figure \ref{F:Q(5,3)}. This slit travels north into the bottommost cylinder from the pole. Taking a slit that travels south into the same cylinder from the zero handles the mergeable component of $\cQ(9,-1)$. Likewise, the analogous construction using a cylindrically stable surface in $\cQ(2,-1,-1)$ handles the case of $\cQ(3,2,-1)$. 

To handle both components of $\cQ(9,3)$, we use that surfaces in the two components can be formed by applying the slit construction to surfaces in the two components of $\cQ(9,-1)$. Every generic cylinder in $\cQ(9,-1)$ is a simple cylinder with the zero of order $9$ appearing on both boundaries (see the results of Masur-Zorich \cite{MZ} as recalled in \cite[Section 4]{ApisaWrightDiamonds}). Use a slit construction with an endpoint of the slit in this simple cylinder to form a simple cylinder on a surface in $\cQ(9,3)$ with the two zeros on opposite boundaries.

% The cases of $\cQ(3,2,-1)$ and $\cQ(9,3)$ can be handled by applying the slit construction to $\cQ(2,-1,-1)$ and $\cQ(9,-1)$. This handles both components of $\cQ(9,3)$ since these components can be reached by applying the slit construction to the two components of $\cQ(9,-1)$. The case of the mergable component of $\cQ(9,-1)$ can be handled by applying the slit construction to $\cQ(5,-1)$. 

To handle $\cQ(5,4,-1)$, we start with a surface in $\cQ(4,1,-1)$ as in Figure \ref{F:Q(4,1,-1)}, then we apply the slit construction using a slit that is contained in the bottommost cylinder and that travels north. The case of $\cQ(5,2,2,-1)$ is handled by splitting the zero $(4)$ into zeros $(2,2)$. As in Kontsevich-Zorich \cite{KZ}, this is a local rel operation that adds in a short saddle connection. When the saddle connection is chosen to have purely imaginary period, the resulting surface remains cylindrically stable. Finally, since the operation only modifies the metric in a neighborhood of the zero of order $4$, it does not affect the cylinder that we use to merge the zero of order $5$ with the pole. 
%\ann{P: Added some detail. It's still a bit hand-wavy, but maybe it's fine.}
\end{proof}

 \begin{figure}
     \centering
     \includegraphics[width =.65\textwidth]{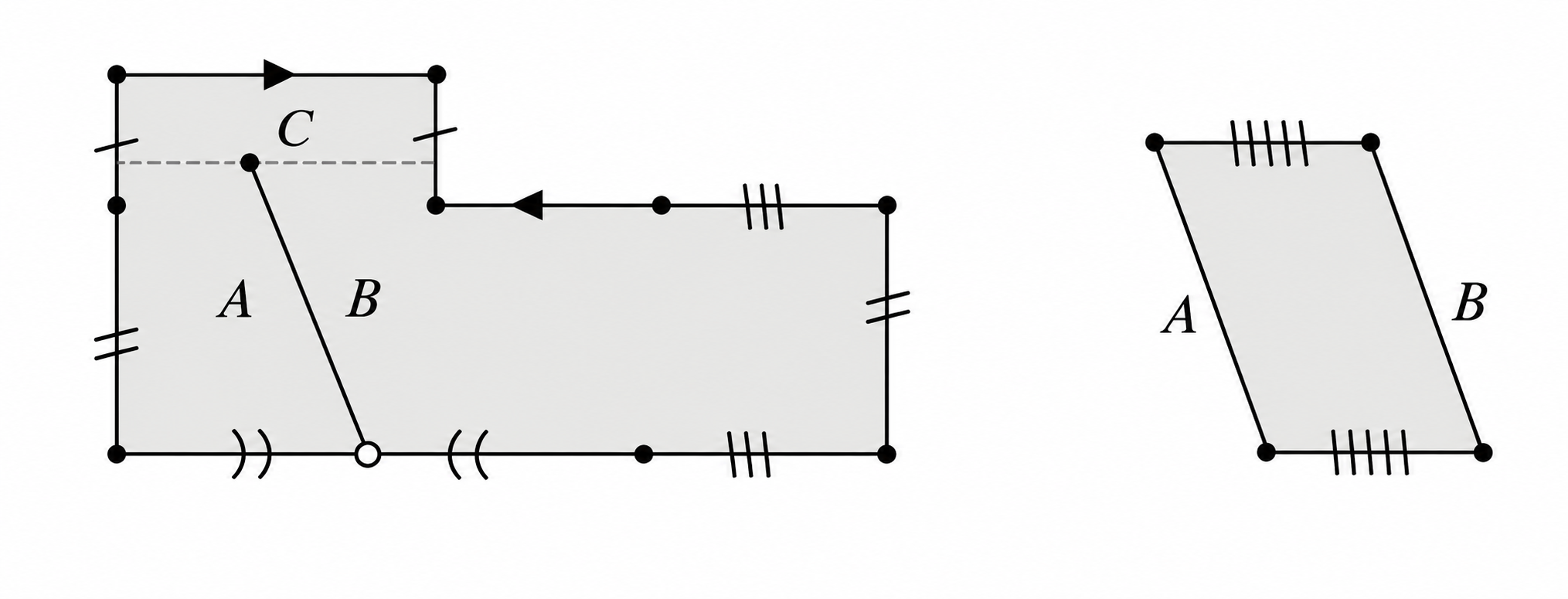}
     \caption{Surface in $\cQ(5,3)$ with a simple cylinder $C$ between the two zeros} 
     \label{F:Q(5,3)}
 \end{figure} 

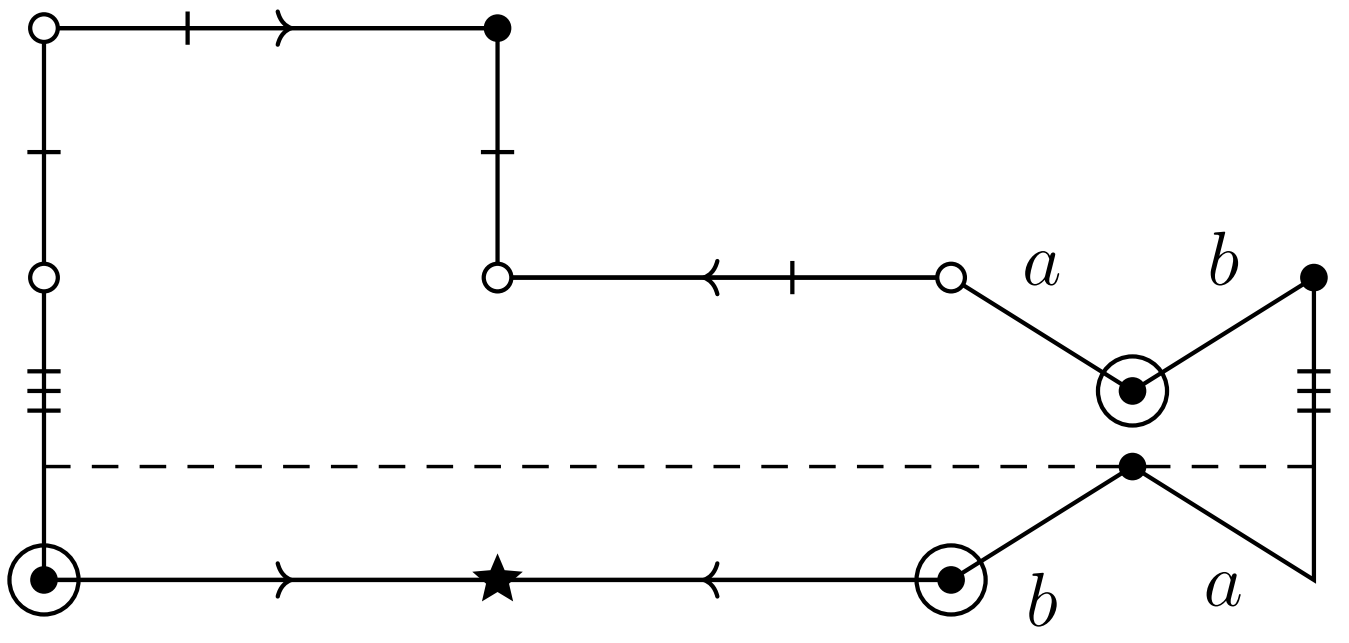
\begin{figure}
     \centering
\begin{tikzpicture}[x=1cm,y=1cm,line width=.8pt,
  every node/.style={font=\large},
  directed/.style={postaction={decorate},decoration={markings,
    mark=at position .55 with {\arrow{>}}}}]
\coordinate (A) at (0,0);
\coordinate (B) at (0,2);
\coordinate (C) at (0,3.65);
\coordinate (D) at (3,3.65);
\coordinate (E) at (3,2);
\coordinate (F) at (6,2);
\coordinate (G) at (7.2,1.25);
\coordinate (H) at (8.4,2);
\coordinate (I) at (8.4,0);
\coordinate (J) at (7.2,.75);
\coordinate (K) at (6,0);
% Polygon boundary, with a white interior suitable for monochrome printing.
\draw (A)--(B)--(C)--(D)--(E)--(F)--(G)--(H)--(I)--(J)--(K)--cycle;
% Dashed interior segment.
\draw[dash pattern=on 5pt off 4pt,line width=.65pt] (0,.75)--(8.4,.75);
% Oriented edge identifications.
\draw[directed] (C)--(D);
\draw[directed] (F)--(E);
\draw[directed] (A)--(3,0);
\draw[directed] (K)--(3,0);
% Single transverse marks; arrows distinguish the horizontal pair.
\draw (.95,3.54)--(.95,3.76);
\draw (4.95,1.89)--(4.95,2.11);
\draw (-.11,2.83)--(.11,2.83);
\draw (2.89,2.83)--(3.11,2.83);
% Triple transverse marks on the two remaining vertical sides.
\foreach \y in {1.12,1.25,1.38}{
  \draw (-.11,\y)--(.11,\y);
  \draw (8.29,\y)--(8.51,\y);
}
% Paired sloping edges. Labels remain upright.
\path (F)--(G) node[midway,above=5pt] {$a$};
\path (G)--(H) node[midway,above=5pt] {$b$};
\path (K)--(J) node[midway,below=5pt] {$b$};
\path (J)--(I) node[midway,below=5pt] {$a$};
% Open and filled vertex markings, as in the sketch.
\foreach \p in {B,C,E,F} \draw[fill=white] (\p) circle (2.6pt);
\foreach \p in {A,D,G,H,J,K} \fill (\p) circle (2.6pt);
\foreach \p in {A,G,K} \draw (\p) circle (6.5pt);
% Starred point on the folded lower boundary.
\node[star,star points=5,star point ratio=2.2,fill=black,
      minimum size=10pt,inner sep=0pt] at (3,0) {};
\end{tikzpicture}
     \caption{The desired surface in $\cQ(4,1,-1)$. The pole is indicated with a star. The circled zero has order $1$. The desired cylinder has one simple boundary, which is the dashed saddle connection.} 
     \label{F:Q(4,1,-1)}
\end{figure}

\appendix

\appendix

\section{Background Material}\label{A:Background}

We begin with a direct proof of Kra's theorem \cite{Kra1981} following McMullen \cite[Theorem 4.1]{Mc}.

\begin{thm}[Kra, McMullen]\label{T:KraMcMullen}
Any quadratic differential in $T^* \cT_{g',n'}$ whose zeros (not occurring at punctures) are of even order generates a Teichm\"uller disk that is a holomorphic retract of Teichm\"uller space.
\end{thm}
\begin{proof}
If $(Y, q)$ is the quadratic differential, then its holonomy cover $(X, \omega)$ is only branched over punctures. So this determines an isometric embedding of $\mathrm{Teich}_{g',n'}$ into $\mathrm{Teich}_{g,n}$ where $g$ is the genus of $X$ and $n$ is the number of preimages of punctures. Fix an isotropic subspace $(a_2, \hdots, a_g)$ of $\ker(\omega) \subseteq H_1(X; \bR)$ and complete it to a Lagrangian subspace by adding in a class $a_1$ where $\int_{a_1} \omega = 1$. Add in classes $(b_1, \hdots, b_g)$ to complete this to a symplectic basis. 

The Teichm\"uller disk of $(X, \omega)$ can be identified with the image of $\bH$ in $\mathrm{Teich}_{g, n}$ given by sending $t := x+iy$ to $X_t$ where $(X_t, \omega_t) := A_t  \cdot (X, \omega)$ where $A_t := \begin{pmatrix} 1 & x \\ 0 & y \end{pmatrix}$. By applying, $A_t$ we may assume that $\int_{b_1} \omega = i$.

Our holomorphic retract $f: \mathrm{Teich}_{g,n} \ra \bH$ will be given by taking a point $X'$ in Teichm\"uller space, finding a basis of holomorphic $1$-forms $\omega_1, \hdots, \omega_g$ so that $\int_{a_i} \omega_j = \delta_{ij}$ and then defining $f(X') = \int_{b_1} \omega_1$. If $\tau_{X'}$ is the period matrix of $X'$ in the indicated symplectic basis, then its imaginary part is positive definite and, since our function returns $(1, 0, \hdots, 0) \tau_{X'} (1, 0, \hdots, 0)^T$ our function is always valued in the upper half plane. Finally, since $\omega_t(a_i) = 0$ if $i > 1$ and $A_t \int_{a_1} \omega = 1$ when $i = 1$, it follows that on $X_t$, $\omega_1 = \omega_t$. So $f(X_t) = \int_{b_1} \omega_t = t$ as desired. 
\end{proof}

%%%%%%%%%%
%
% DO NOT DELETE:
% Explicit understanding of the holomorphic retraction!
%
%%%%%%%%
%
% Although this won't be needed suppose that $(X, \omega)$ in our previous proof was horizontally periodic with horizontal core curves spanning a Lagrangian subspace of homology. Let $(a_1, \hdots, a_g)$ be a basis of this subspace. In the Teichm\"uller polyplane, the periods of these curves is constant and, after rescaling, we may assume that the period of $a_1$ is $1$. By changing basis, we may assume that $\int_{a_i} \omega = 0$ when $i > 1$ for any $\omega$ in the polyplane. Let $(b_1, \hdots, b_g)$ be the completion of $(a_1, \hdots, a_g)$ to a symplectic basis. Then our retraction on the polyplane is simply given by sending $(X', \omega')$ in the polyplane to $\int_{b_1} \omega'$. 

We will now prove the key ingredient in Markovic's criterion, i.e. Theorem \ref{thm:retraction on polyplane}.

\begin{thm}\cite[Theorem 4.1]{Gekhtman-Polyplane}\label{L:ConvergenceToLinear}
Fix $a_1, \hdots, a_k \in \bR_{>0}$ that sum to $1$. Let $\cC$ be the space of holomorphic functions $f: \bH^k \ra \bH$ so that $f(t,\hdots, t) = t$ and $(\del_i f)(t, \hdots, t) = a_i$. Note that, with the compact-open topology, $\cC$ is compact\footnote{Compactness follows by Arzela-Ascoli and the fact that these maps are contractions in the Kobayashi metric.} and convex. Define $f_t := f(z_1-t, \hdots, z_k-t) + t$, which is a flow on $\cC$. Set $g := \sum_i a_i z_i$. Let $U$ be a neighborhood of $g$, then, for any $f \in \cC$, $(f_t)_{t \geq 0}$ spends asymptotically all of its time in $U$.
\end{thm}
\begin{proof}
Let $\mu_r$ be the pushforward of normalized Lebesgue on $[-r, r]$ to $(f_t)_{t \in [-r, r]}$. Let $x_r$ be the convex linear combination of points in $\cC$ using $\mu_r$, i.e. the barycenter of $\mu_r$. 

We will show that $x_r$ converges to $g$ and that $g$ is an extreme point of $\cC$, but let's see first why this is sufficient. Since barycenters of weakly converging measures converge, pass to a weakly converging subsequence $\mu_{r_n} \ra \mu$. Since the barycenter of $\mu$ is an extreme point, i.e. $g$, it follows that $\mu = \delta_g$. So $\mu_r \ra \delta_g$ since all weakly convergent subsequences converge to the same limit. The claim about the orbit spending asymptotically all of its time in $U$ now follows.

\begin{lem}
$g$ is an extreme point.
\end{lem}
\begin{proof}
The generalized Schwarz lemma says that for a map $f: \bH^n \ra \bH$ and for any $z \in \bH^n$, $\sum_i \mathrm{Im}(z_i)|\del_i f(z)| \leq \mathrm{Im}(f(z))$. When $f$ is a real linear combination of $z_i$, equality holds. So if $g = \lambda f + (1-\lambda)h$ then
\[ \mathrm{Im}(g(z)) = \sum_i \mathrm{Im}(z_i)|\lambda \del_i f(z) + (1-\lambda)\del_i h(z)| \leq \sum_i  \mathrm{Im}(z_i) \left( \lambda|\del_i f(z)| + (1-\lambda)|\del_i h(z)| \right)).  \]
Continuing,
\[ \mathrm{Im}(g(z)) \leq \lambda \mathrm{Im}(f(z)) + (1-\lambda)\mathrm{Im}(h(z)) = \mathrm{Im}(g(z)).  \]
Since each inequality is an equality, $\del_i f(z)$ and $\del_i h(z)$ (and hence $\del_i g(z))$ have the same argument for any $z$. So their ratio is real-valued and hence constant by the open mapping theorem. This constant is $1$ by considering partials on the diagonal. So $f = h = g$.
\end{proof}

Since any convergent subsequence of $x_r$ as $r \ra \infty$ converges to a flow invariant point in $\cC$, it suffices to show the following.

\begin{lem}\label{L:FlowInvariant}
$g$ is the only flow-invariant function in $\cC$.
\end{lem}
\begin{proof}
Let $f$ be a flow-invariant function in $\cC$ and consider $h := f - g$. For any $w = (w_1, \hdots, w_k) \in \bH^k$ and any $z \in \bC$ notice that $h(w_1 - z, \hdots, w_k - z) = h(w)$ whenever $(w_1 - z, \hdots, w_k - z) \in \bH^k$. Now define a function $H: \bC^{k-1} \ra \bC$ as follows. Given $w' = (w_1, \hdots, w_{k-1})$ choose $\tau \in \bH$ so that $w_i + \tau \in \bH$ for all $i$. Then set $H(w') := h(w_1 + \tau, \hdots, w_{k-1} + \tau, \tau)$. This function is entire. Since $g+h$ sends $\bH^n$ to $\bH$ it follows that
\[ \mathrm{Im}(H(w_1, \hdots, w_{k-1})) \geq a_k\tau + \sum_i a_i (\mathrm{Im}(w_i)+\tau) \geq 2\tau \]
where $\tau = -\min_j(\mathrm{Im}(w_j)) = \| \mathrm{Im}(-w) \|_\infty$. Rephrased, $-\frac{1}{2} H(w)$ is an entire function whose imaginary part is bounded above by $\|w\|_\infty$. So by Liouville's theorem, $H$ is affine. By assumption $H(0, \hdots, 0) = 0$, so $H$ is linear. So $f$ is linear and hence completely specified by its partials, so $f = g$.
\end{proof}
\end{proof}

\begin{proof}[Proof of Theorem \ref{thm:retraction on polyplane}:] Let $(X, q)$ be the unit-area horizontally periodic quadratic differential whose Teichm\"uller disk is a holomorphic retract. Let $F: \mathrm{Teich}_{g,m} \ra \mathbb{H}$ be the retract and let $f: \bH^k \ra \bH$ be its restriction to the polyplane. We identify $(X, q)$ with $(i, \hdots, i)$ and suppose that $f(i,  \hdots, i) = i$. By classical Teichm\"uller theory, $f$ satisfies the assumptions of Theorem \ref{L:ConvergenceToLinear} where $a_i$ is the area of the $i$th cylinder on $X$ (see \cite[Lemma 4.2]{GM2020}). Note that $f_t = t + f \circ u_{-t}$ where $u_t$ is the matrix that applies $\begin{pmatrix} 1 & -t/h_i \\ 0 & 1 \end{pmatrix}$ to the $i$th cylinder, where $h_i$ is the height of that cylinder. In the stratum, the $u_t$ flow is a linear flow on a torus and so, for any positive integer $n$, it spends a definite fraction of time within $1/n$ of its starting point. Let $(U_n)$ be a neighborhood basis of $g := \sum_i a_i z_i$ in $\cC$, i.e. a decreasing sequence of open sets so that $\bigcap_n U_n = \{g\}$. By Theorem \ref{L:ConvergenceToLinear}, we may find a time $t_n$ so that $f_{t_n}$ belongs to $U_n$ and so $u_{t_n}$ is $\frac{1}{n}$ close (in the polyplane) to a Dehn multitwist $\tau_n := T_1^{m_{1,n}} \hdots T_k^{m_{k,n}}$ where $T_i$ is a Dehn twist in the $i$th cylinder. By construction, $f_{t_n}$ converges to $g$. So $t_n + F \circ \tau_n$ is a holomorphic map from Teichm\"uller space to $\bH$ whose restriction to the polyplane is $t_n + f(z_1 - \frac{m_{1,n} \ell_1}{h_1}, \hdots, z_k - \frac{m_{k,n} \ell_k}{h_k})$ where $\ell_i$ is the circumference of the $i$th cylinder. By assumption $m_{i,n}\ell_i/h_i$ is $\frac{1}{n}$-close to $t_n/h_i$ for all $i$. This shows that the restriction of $t_n + F \circ \tau_n$ is locally uniformly asymptotic to $f_{t_n}$ (using that holomorphic maps are retractions for the Kobayashi metric). By Arzela-Ascoli, it follows that $t_n + F \circ \tau_n$ has a convergent subsequence whose restriction to the polyplane converges to $g$.
\end{proof}

\bibliography{mybib}{}
\bibliographystyle{amsalpha}
\end{document}